\documentclass[12pt,thmsa]{amsart}
\usepackage{amssymb}

\def\Cal{\mathcal}

\def\B{{\Cal B}}

\def\P{{\Cal P}}

\def\S{{\Cal S}}
\def\F{{\Cal F}}
\def\K{{\Cal K}}
\def\I{{\Cal I}}

\def\T{{\Cal T}}
\def\W{{\Cal W}}
\def\J{{\Cal J}}

\def\d{\partial}
\def\tr{{\hbox{\rm tr}}}

\def\Ma{\frM_{n,m}}
\def\Mt{\frM_{n-k,m}}
\def\Mkm{\frM_{k,m}}

\def\Mmk{\frM_{m,k}}

\def\W{\mathcal{W}}

\def\Z{\mathcal{Z}}

\def\f0{f_0}
\def\Fc0{\varphi_0}
\def\rn{\bbr^n}

\def\I_k {I_{-}^{k/2}}
\def\I+k {I_{+}^{k/2}}

\def\vnk{V_{n,n-k}}
\def\vnm{V_{n,m}}
\def\cd{\vnk\times\frM_{n-k,m}}

\def\sigk{\sig_{n,n-k}}

\def\bbr{{\Bbb R}}

\def\bbn{{\Bbb N}}

\def\bbc{{\Bbb C}}

\def\rank{{\hbox{\rm rank}}}
\def\diag{{\hbox{\rm diag}}}

\def\supp{{\hbox{\rm supp}}}

\def\tr{{\hbox{\rm tr}}}

\def\const{{\hbox{\rm const}}}

\def\det{{\hbox{\rm det}}}

\def\Pr{{\hbox{\rm Pr}}}

\def\fc{\varphi(\xi,t)}

\def\part{\partial}
\def\intl{\int\limits}
\def\b{\beta}

\def\Gam{\Gamma}

\def\a{\alpha}
\def\om{\omega}

\def\Del{\Delta}
\def\del{\delta}
\def\vp{\varphi}

\def\g{\gamma}
\def\gam{\gamma}

\def\sig{\sigma}
\def\lam{\lambda}
\def\z{\zeta}
\def\e{\varepsilon}
\def\t{\tau}
\def\eq{\xi 'x=t}

\def\rf{\hat f (\xi,t)}
\def\df{\check \varphi (x)}

\font\frak=eufm10

\def\fr#1{\hbox{\frak #1}}

\def\frM{\fr{M}}

\def\frT{\fr{T}}

\def\pk{\P_k}

\def\vmk{V_{m,k}}
\def\gma{\Gamma_m(\a)}

\def\gk{\Gamma_k}

\def\const{{\hbox{\rm const}}}

\def\det{{\hbox{\rm det}}}

\def\p{\P_m}
\def\gm{\Gamma_m}
\def\tr{{\hbox{\rm tr}}}
\def\part{\partial}
\def\intl{\int\limits}
\def\b{\beta}

\def\Gam{\Gamma}

\def\a{\alpha}
\def\cpm{\overline\P_m}

\newtheorem{theorem}{Theorem}[section]
\newtheorem{lemma}[theorem]{Lemma}

\theoremstyle{definition}
\newtheorem{definition}[theorem]{Definition}

\newtheorem{example}[theorem]{Example}

\theoremstyle{corollary}
\newtheorem{corollary}[theorem]{Corollary}
\theoremstyle{remark}
\newtheorem{remark}[theorem]{Remark}

\numberwithin{equation}{section}

\newcommand{\be}{\begin{equation}}
\newcommand{\ee}{\end{equation}}

\newcommand{\bea}{\begin{eqnarray}}
\newcommand{\eea}{\end{eqnarray}}
\newcommand{\Bea}{\begin{eqnarray*}}
\newcommand{\Eea}{\end{eqnarray*}}

\begin{document}

\title [ Zeta integrals ]{Zeta integrals and integral geometry  in the space of rectangular matrices}

\author{Boris Rubin}



\address{Institute of Mathematics,
 Hebrew University, Jerusalem 91904, \newline Israel}

\email{boris@math.huji.ac.il}

\thanks{ The work  was supported in part by
the Edmund Landau Center for Research in Mathematical Analysis
and Related Areas, sponsored by the Minerva Foundation (Germany).}





\begin{abstract}
The paper is devoted to interrelation between the  zeta distribution
 $\Z(f,\a-n)= \int f(x) \,
\det(x'x)^{(\a -n)/2}dx$  and the Radon transform  on the space
$\Ma$  of $n \times m$ real matrices  $x=(x_{i,j})$.  We present a
self-contained proof of the Fourier transform formula for  this
distribution. Our method differs from that of J. Faraut and A.
Kor\'anyi [FK] in the part related to justification of the
corresponding Bernstein identity. We suggest a new proof of this
identity based on explicit representation of the radial part of
the Cayley-Laplace operator $\Del=\det(\partial'\partial), \;
\partial= (\partial_{i,j})_{n \times m}$.
 We also study convolutions with  normalized  zeta distributions,
 and the corresponding  Riesz potentials. The results are applied  to investigation of Radon transforms  on  $\Ma$.
 \end{abstract}

\maketitle

\centerline{Contents} \centerline{} 1. Introduction.

2. Preliminaries.

3. Radial functions and the Cayley-Laplace operator.

4. Zeta integrals.

5. Convolutions  with zeta distributions and Riesz
potentials.

6. Radon transforms.

7. Appendix.

\section{Introduction}

\setcounter{equation}{0}

 Diverse problems of integral geometry in the $N$-dimensional real space
$\bbr^N$  get a new  flavor if the dimension $N$ has the form
 $N=nm$, and $\bbr^N$ is identified with the space
  $\Ma$  of $n \times m$
real matrices  $x=(x_{i,j})$. If we accept this point of view, a
number of
 new ``higher rank'' phenomena will come into play. In the present
 paper we focus on analytic tools which allow to investigate the properly defined
  Radon transform $f(x) \to \hat f(\t)$ on  $\Ma$. This transform assigns
to a function $f(x)$ on  $\Ma$ a collection of integrals of $f$
over the so-called {\it matrix $k$-planes} $\t$ in   $\Ma$.
 Each such plane is, in fact, an
ordinary $km$-dimensional plane in $\bbr^{nm}$, but the manifold  $\frT$
 of all
 matrix $k$-planes is essentially ``smaller'' than the manifold of all
 $km$-dimensional planes in $\bbr^{nm}$.  The main problem is to
 reconstruct $f(x)$ from  known data  $\hat f(\t)$, $\t \in \frT$, for possibly large class of functions $f$.

In the rank-one case $m=1$ this problem is well investigated; see,
e.g., [E], [Hel], [Ru2], and references therein. The standard tools are
the  Fourier analysis, convolution operators, and the
corresponding real-variable techniques, including approximation to
the identity, singular integral operators, and fractional
integrals. In the higher rank case $m \ge 2$, all these tools
still require
 a further development.

Let us pass to details. Suppose that
  $f(x)$ is a Schwartz function on $\frM_{n,m}$, $n\geq m$.  
The Fourier transform
 of   $f$ is
defined by \be (\F f)(y)=\intl_{\Ma} \exp(\tr(iy'x)) f (x)
dx,\qquad y=(y_{i,j})\in\Ma,\ee
where   $y'$ denotes the transpose of
 $y$,
$dx=\prod^{n}_{i=1}\prod^{m}_{j=1}
 dx_{i,j}$. This is the ordinary Fourier transform on $\bbr^{nm}$. Given a complex number $\a$, we denote \be\label{zet}
\Z(f,\a-n)=\intl_{\Ma} f(x) |x|^{\a-n}_m dx \ee where $|x|_m =\det
(x'x)^{1/2}$ is the volume of the parallelepiped spanned by the
column-vectors of the matrix $x$. For $m=1$, $|x|_m$ is the usual euclidean norm on $\bbr^n$.
 Integrals like (\ref{zet})
 are known as  {\it  zeta integrals}, the $\S'$-distribution defined
 by analytic continuation of (\ref{zet}) is  called  {\it a  zeta
distribution}, and
the corresponding
 function $\a \to\Z(f,\a-n) $ is sometimes  called  {\it a zeta function}.

We begin our study by giving a relatively simple self-contained
proof of the well-known functional relation
 \be\label{z-f11}
\frac{\Z(f, \a-n)}{\gm(\a/2)}=\pi^{-nm/2} 2^{m(\a-n)}\frac{\Z(\F
f, -\a)}{\gm((n-\a)/2)}, \ee where the normalizing denominators
are Siegel gamma functions associated to the cone  of positive
definite $m \times m$ matrices. Both sides of (\ref{z-f11}) are
understood in the sense of analytic continuation.
 The equality (\ref{z-f11}) is of fundamental importance for the sequel. It 
 is, in fact,  a Fourier transform formula for
the zeta distribution. Convolutions with zeta distributions
coincide (up to normalization) with Riesz potentials on $\Ma$
playing a vital role in integral geometry and analysis on matrix spaces [OR1],
[P], [Sh].

Our interest to this topic  is also motivated by the following.
 A number of  proofs of the equality (\ref{z-f11})
and its generalizations can be found in the literature; see, e.g., [Cl],
 [Far], [FK], [Ge], [Ra], [Sh]. Unfortunately, these proofs are not ``self-contained enough'' and
 some important technical details  are skipped.
These details are crucial by taking into account an essential
difference between
 the case $2m<n+2$ and   $2m\geq n+2$. This difference is not
 indicated in some papers. The crux  is that the distribution on the left-hand side of (\ref{z-f11})
is regular if and only if $Re\, \a >m-1$, whereas the right-hand
side is regular if and only if $Re\, \a <n-m+1$. For $2m\geq n+2$
these two sets on the complex plane are separated one from another!
 In the  case $2m<n+2$ the
 proof is elementary [OR1]. For $2m\geq n+2$,  justification  of
 (\ref{z-f11}) represents a difficult problem.

 The phenomenon of
   lack of common domain of
 regularity was investigated by E. Stein [St1] for Riesz distributions on the space of complex $n \times n$
 matrices. For $n=m$, the formula (\ref{z-f11})  includes the result of S. S.
Gelbart \cite[Chapter IV] {Ge} and a special case of Proposition
II-9 of M.  Ra\"{\i}s \cite{Ra}.  For $n\geq m$, this formula can
be found in \cite{Far}, \cite[Theorem XVI, 4.3]{FK}, and \cite{Cl}
in a more general context of representations of Jordan algebras;
see also [Sh,  p. 108]. For more information on zeta
distributions, see [Ig], [SS], [Shin], and references therein.

 Our proof of (\ref{z-f11}) utilizes the idea from
\cite{FK} to employ $\K$-Bessel functions and the so-called
Bernstein identity \be \label{ber}\Del |x|_m^{\lam}=\B(\lam)
|x|_m^{\lam-2},\ee
 \[ \B(\lam)=(-1)^m
\prod\limits_{i=0}^{m-1}(\lam+i)(2-n-\lam+i), \]
 where $\Del$ is  the
Cayley-Laplace operator defined by $\Del=\det(\partial'\partial), \;
\partial= (\partial_{i,j})_{n \times m}$. This identity amounts to pioneering papers of J. Bernstein; see, e.g., [B], [Ig].
  A  novelty of  our approach  is
that we first derive an explicit formula for the radial part of
$\Del$, and then use it for justification  of (\ref{ber}). It is
worth noting that the proof of the Bernstein identity in
\cite{FK} and \cite{Cl} employs the heat kernel and integration by
parts over the corresponding cone. Justification of this
integration by parts
 was unfortunately skipped. That was one of the reasons why we tried to find
an alternative proof.
 We hope that our method  can be
extended to more general settings similar to those in  \cite{FK} and \cite{Cl}.

The second concern of the paper is  the so-called {\it Wallach set} $\W$ of
 the normalized
 zeta distribution \be\label{nnzd} (\z_\a, f)=a.c.
\frac{1}{\gm(\a/2)}\intl_{\frM_{n,m}} \! f(x)|x|_m^{\a -n} dx, \qquad
f \in \S(\Ma),\ee
where``$a.c.$'' abbreviates analytic continuation. This is
  an entire function of $\a$. The set $\W$ consists of two
  components, one of which is continuous and another is discrete.
   The continuous component $\{\a: Re \, \a > m-1 \}$ is
  constituted by those $\a \in \bbc$ for which the distribution (\ref{nnzd}) is regular, that is,
the integral in (\ref{nnzd})  absolutely converges. The discrete
 component  consists of non-negative integers $\a=k,\; k=0,1,
\ldots , m-1 $ (outside of the domain of regularity), for which
the distribution (\ref{nnzd}) is a positive measure. We obtain
explicit representation  of  $\z_\a$ for all $\a \in \W$, and
investigate convolutions $(\z_\a * f)(x)$ assuming
 $f \in L^p(\Ma)$. In particular, under some natural restrictions on  $n$ and $m$,
we  prove that these convolutions are well defined almost
everywhere on $\Ma$ provided \be\label{st} 1\le p<\frac{n}{Re \,
\a +m-1}. \ee This result agrees with the case $m=1$ in [St2]. The
question whether the condition (\ref{st}) is sharp remains open.

Another normalization of the integral (\ref{zet}) gives rise to
 the Riesz potential  \be\label{rss} (I^\a
f)(x)=\frac{1}{\gam_{n,m} (\a)} \intl_{\Ma} f(x-y) |y|^{\a-n}_m
dy, \ee
\[ \gam_{n,m} (\a)=\frac{2^{\a m} \, \pi^{nm/2}\, \Gam_m (\a/2)}{\Gam_m
((n-\a)/2)},   \quad \a \neq n-m+1, \,
n-m+2, \ldots \, ,\]
the Fourier transform of which is $|y|_m^{-\a}(\F f)(y)$ in a suitable sense.
 The operator $I^\a$ can be regarded as the $(-\a/2)$th power of
 $(-1)^m\Del$ where $\Del$ is the Cayley-Laplace
operator.

Our investigation of  convolutions $(\z_\a * f)(x)$ and the  Riesz
potentials $I^\a f$ employs  matrix modification of the
Gauss-Weierstrass integral (see [SW], [Ta] for $m=1$) defined by
\be\label{gaw} (W_t f)(x)=\intl_{\Ma} h_t (x-y) f(y) dy, \ee where
$t$ is a positive definite $m \times m$ matrix, and

\be\label{hit} h_t (x)=(4\pi)^{-nm/2} \, \det(t)^{-n/2} \exp (-\tr
(t^{-1} x'x)/4) \ee is
 the relevant heat kernel. This approach allows to avoid
 essential technical difficulties. The idea is  to
represent the Riesz potential through  the   lower-dimensional
 G{\aa}rding-Gindikin fractional integral  by a simple formula
\be\label{rggg}  (W_t [I^\a f])(x) =(I_{-}^{\a/2}g_x)(t), \qquad
g_x(t) =(W_{t} f)(x), \ee provided $ m-1<Re \, \a<n-m+1$. We
recall that the G{\aa}rding-Gindikin fractional integral
 on the cone $\p$ of positive
 definite matrices has the form
  \be\label{gig} (I_{-}^\lam g)(t) = \frac {1}{\Gam_m (\lam)}
\intl_t^\infty g(\t) \, \det(\tau -t)^{\lam-d} d\t, \ee where
$d=(m+1)/2, \; Re \, \lam > d-1$, and integration is performed
over all $\t \in \p$ so that $\t -t \in \p$ . These integrals were
introduced by Lars G{\aa}rding \cite{Ga1} who wrote  \cite{Ga2}:

``... Actually the origin of my integral was a statistics paper
(samples of mean, variance etc. from a multivariate Gaussian
distribution). The idea to use it in analysis came from my many
meetings with my beloved teacher Marcel Riesz".

Integrals (\ref{gig}) were substantially  generalized by S.
Gindikin \cite{Gi}. They have proved to be an important tool in
PDE \cite{VG}, \cite {Rab}, and in the theory of Radon transforms
on Grassmann manifolds and matrix spaces \cite{GR}, [OR2],
\cite{Ru3}. If $t=0$, then $W_t$ turns into the identity operator
and (\ref{rggg}) reads \be\label{rrg} (I^\a f)(x) = \frac
{1}{\Gam_m(\a/2)} \intl_{\p} |t|^{\a/2-d}(W_t f)(x) dt.\ee

Our application of zeta distributions to integral geometry is
based on intimate connection between the Riesz potentials
(\ref{rss}) and the matrix $k$-plane Radon transform $f \to \hat
f$ on $\Ma$. This connection is realized via the generalized
Fuglede formula \be\label{fff}  (\hat f)^\vee (x)=\const \, (I^k
f)(x)\ee which is well known for $m=1$ [Fu], [Hel]. Here the
left-hand side is the mean value of $\hat f(\t)$ over all  matrix
$k$-planes $\t$ ``passing through $x$''. For sufficiently good
functions $f$, this formula was established in [OR1]. We  justify
it for all   $f \in L^p(\Ma)$, $1 \le p<n/(k+m-1)$, and  apply to
the inversion problem for the Radon transform
 $f \to \hat f$. Solution to this problem is presented in the framework
 of the relevant theory of distributions. For $m=1$, such a theory was developed by  V.I. Semyanistyi [Se].

The paper is organized as follows. In Section 2 we establish our
notation and recall some basic facts important for the sequel. In
Section 3 we derive a formula for the radial part of the
Cayley-Laplace operator and prove the Bernstein identity
(\ref{ber}). In Section 4 we
 study the zeta integral (\ref{zet}), prove (\ref{z-f11}), and obtain
explicit
representations of the  normalized
 zeta distribution (\ref{nnzd})
 in the discrete part of the Wallach set.  Convolutions  with zeta distributions and Riesz
potentials of $L^p$ functions are investigated in Section 5.
Section 6 is devoted to applications of the results of preceding
sections to integral geometry. We define the matrix $k$-plane
Radon transform and justify the generalized Fuglede formula
(\ref{fff}) for  $f \in L^p(\Ma)$, $1 \le p<n/(k+m-1)$. Then we
introduce the space $\Phi'$ of distributions of the Semyanistyi
type. Two $\Phi'$-distributions coincide if and only if their
Fourier transforms differ by a tempered distribution supported by
the singular set $\{y: y \in \Ma, \rank (y)<m \}$. An inversion
formula for the Radon transform $f \to \hat f$ for $f\in L^p(\Ma)$
is then obtained in the $\Phi'$-sense. For convenience of the
reader, evaluation of some useful integrals is presented in
Appendix.

The paper contains a number of open problems that might be of
interest. These are stated in Remarks
\ref{hl}, \ref{rpo}, and \ref{op3}.

{\bf Acknowledgements.} I am grateful to Prof. Jacques Faraut for
valuable remarks, and Prof. Jean-Louis Clerc for correspondence.
Special thanks go to Prof. Evgenyi Evgenyievich Petrov whose
pioneering paper [P] and  friendly comments  were so
inspiring. I am also indebted to Dr. Sergei Khekalo and Dr. Elena
Ournycheva for useful discussions.

{\section{Preliminaries}

\setcounter{equation}{0}

\subsection{Matrix spaces}\label{s2.1}

  Let $\frM_{n,m}$ be the
space of real matrices having $n$ rows and $m$
 columns.  We identify $\frM_{n,m}$
 with the real Euclidean space $\bbr^{nm}$.
The letters $x$, $y$, $r$, $s$, etc. stand for both the matrices
and points since it is always clear from the context which is
meant. If $x=(x_{i,j})\in\Ma$, we write
$dx=\prod^{n}_{i=1}\prod^{m}_{j=1}
 dx_{i,j}$ for the elementary volume in $\Ma$. In the following
  $x'$ denotes the transpose of  $x$, $I_m$ is the identity $m \times m$
  matrix,  $0$ stands for zero entries. Given a square matrix $a$,  we denote by $|a|$ the absolute value of the determinant of $a$;  $\tr (a)$
  stands for  the trace  of $a$.

Let $\S_m$ be the space of $m \times m$ real symmetric matrices
$s=(s_{i,j}),
 \, s_{i,j}=s_{j,i}$. It is a measure space   isomorphic to $\bbr^{m(m+1)/2}$
 with the volume element $ds=\prod_{i \le j} ds_{i,j}$.
We denote by  $\p$ the open  convex cone of positive definite
matrices in $\S_m$; $\bar \P_m$ is the closure of $\p$ that
consists of positive semi-definite matrices.
  For $r\in\p$  we
write $r>0$. The inequality $r_1 > r_2$  means $r_1 - r_2 \in
 \p$. If $a$ and $b$ are  positive semi-definite matrices, the
symbol $\int_a^b f(s) ds$ denotes
 integration over the set
$$
 \{s : s \in \p, \, a<s<b \}=\{s : s-a \in \p, \, b-s \in \p \}.
$$ The group $G=GL(m,\bbr)$ of
 real non-singular $m \times m$ matrices $g$ acts on $\p$ transitively
  by the rule $r \to g'rg$.  The corresponding $G$-invariant
 measure on $\p$ is
\be\label{2.1}
  d_{*} r = |r|^{-d} dr, \qquad |r|=\det (r), \qquad d= (m+1)/2 \ee \cite[p.
  18]{Te}. The cone $\p$ is a $G$-orbit in $\S_m$ of the identity matrix $I_m$.
  The boundary $\partial\p$ of $\p$
  is a union   of $G$-orbits of $m\times m$ matrices
\[e_k=\left[\begin{array}{ll}  I_k&0\\
0&0
\end{array}\right], \qquad k=0,1, \ldots , m-1.\]
   More information about the boundary structure of $\p$ can be
found in \cite[p. 72]{FK},  and \cite[ p. 78]{Bar}.

We denote by $T_m$  a group of
 upper triangular matrices
\be\label{2.17}
 t=\left[\begin{array}{ccccc} t_{1,1} & {}   & {}  & {}  & {} \\
                              {} & {.}  & {}  & t_{*}  & {} \\
                              {} & {}   & {.} & {}  & {} \\
                              {} & {0}   & {}  & {.} & {} \\
                               {} & {}   & {}  & {}  & t_{m,m}

\end{array} \right], \qquad t_{i,i} >0, \ee $$t_{*}=\{ t_{i,j}: i<j \} \in \bbr^{m(m-1)/2}.
$$ Each $r \in \p$ has a unique representation $r=t't, \; t \in
T_m$, so that
 \bea \label{2.2}
 \intl_{\p} f(r) dr &=& \intl_0^\infty  t_{1,1}^m \,
dt_{1,1}\intl_0^\infty  t_{2,2}^{m-1} \, dt_{2,2} \,
 \ldots \nonumber \\ &\times & \intl_0^\infty  t_{m,m}  \tilde f (t_{1,1}, \ldots ,t_{m,m}) \,
 dt_{m,m},
\eea
\[
 \tilde{ f} (t_{1,1}, \ldots ,t_{m,m})=2^m  \intl_{\bbr^{m(m-1)/2}} f(t't) \, dt_{*}, \quad dt_{*}=\prod_{i<j}
 dt_{i,j},
\]
 \cite[p. 22]{Te} , \cite[p. 592]{Mu}. In the last integration, the diagonal
 entries of the
 matrix $t$ are given by the arguments of $\tilde f$, and the
 strictly upper
 triangular entries of $t$ are variables of integration.

 Let us recall some useful formulas for
Jacobians.
 \begin{lemma}\label{12.2} {\rm (see, e.g., \cite[ pp. 57--59]
 {Mu})}.\hskip10truecm

\noindent
 {\rm ($i$)} \ If $ \; x=ayb$, where $y\in\Ma, \; a\in  GL(n,\bbr), \; b \in  GL(m,\bbr)$, then
 $dx=|a|^m |b|^ndy$. \\
 {\rm ($ii$)} \ If $ \; r=q'sq$, where $s\in S_m, \; q\in  GL(m,\bbr)$,
  then $dr=|q|^{m+1}ds$. \\
  {\rm ($iii$)} \ If $ \; r=s^{-1}$, where $s\in \p$,   then $r\in
  \p$,
  and $dr=|s|^{-m-1}ds$.
\end{lemma}

  In the following  $(\lam,m)=\lam (\lam +1) \ldots (\lam +m-1)$ is the Pochhammer symbol, $\del_{ij}$ is
 the  Kronecker delta; ``$a.c.$" abbreviates analytic
 continuation. To facilitate presentation, we
shall use the symbols ``${}\simeq {}$" and ``${}\lesssim{}$" ,
instead of ``$=$" and ``$\le$", respectively, to indicate that the
corresponding relation holds up to a constant multiple.
 All standard spaces of functions of matrix argument, say, $x=(x_{i,j})$,  are
 identified with the corresponding spaces of functions of $nm$ variables
$x_{1,1}, x_{1,2}, \ldots, x_{n,m}$. The Fourier transform of a
function  $f\in L^1(\Ma)$ is defined by \be\label{ft} (\F
f)(y)=\intl_{\Ma} \exp(\tr(iy'x)) f (x) dx,\qquad y\in\Ma .\ee

\subsection{Gamma  functions, beta functions, and Bernstein polynomials}\label{s2.2}
  The Siegel gamma
 function associated to the cone $\p$
 is defined by
\be\label{2.4}
 \gm (\a)=\intl_{\p} \exp(-\tr (r)) |r|^{\a -d} dr, \qquad d=(m+1)/2. \ee
It is easy to check \cite[ p. 62]{Mu}, that the integral
(\ref{2.4}) converges absolutely
 if and only if $Re \, \a>d-1$, and represents a product of
 ordinary gamma functions:
\be\label{2.5}
 \gm (\a)=\pi^{m(m-1)/4}\prod\limits_{j=0}^{m-1} \Gam (\a- j/2)\;
 .
 \ee
This implies
 \be\label{Poh} (-1)^m
\frac{\gm(1-\a/2)}{\gm(-\a/2)}=2^{-m}\frac{\Gam(\a+m)}{\Gam(\a)}=2^{-m}(\a,m),
\ee $(\a,m)=\a(\a+1)\cdots (\a+m-1)$ being the Pochhammer symbol.
If $1  \le k<m, \; k \in \bbn, $ then  \be\label{2.5.1}
\gma=\pi^{k(m-k)/2} \Gamma_{k}(\a)\Gamma_{m-k}(\a-k/2), \ee
\be\label{2.5.2}
\frac{\gm(\a)}{\gm(\a+k/2)}=\frac{\gk(\a+(k-m)/2)}{\gk(\a+k/2)}.\ee

The beta function of the cone  $\p$ is defined by

\be\label{2.6}
 B_m (\a ,\b)=\intl_0^{I_m} |r|^{\a -d} |I_m-r|^{\beta -d} dr, \qquad d=(m+1)/2.\ee
This integral converges absolutely if and only if $Re
 \, \a, Re \, \b >d-1$, and obeys the classical relation \cite[ p. 130]{FK}
\be\label{2.6.1}
 B_m (\a ,\b)=\frac{\gm (\a) \gm (\b)}{\gm (\a+\b)}.
\ee

 Let $r = (r_{i,j}) \in  \p$. We define the
 following differential operator acting in the $r$-variable:
 \be\label{2.50} D \equiv D_{ r}=
 \det \left ( \eta_{i,j} \, \frac {\partial}{\partial
 r_{i,j}}\right ), \quad
  \eta_{i,j}= \left\{
 \begin{array} {ll} 1  & \mbox{if $i=j$}\\
1/2 & \mbox{if $i \neq j.$}
\end{array}
\right. \ee

\begin{lemma} {}\hfil

\noindent {\rm (i)}   For  $r \in \P_m$ and $z  \in
\S_m^{\bbc}$ (the complexification of $\S_m$), \be\label{3.41}
D_{r} [\exp(-\tr(rz))]=(-1)^m \det(z)\exp(-\tr(rz)). \ee

\noindent {\rm (ii)}  For  $r \in \P_m$ and $\a  \in \bbc$,
\be\label{3.40} D(|r|^{\a-d}/\gma)= |r|^{\a-1-d}/\gm(\a-1),\qquad
d=(m+1)/2.\ee
\end{lemma}

The proof of this statement is simple and can be found in \cite[p.
813]{Ga1}, \cite[p. 125]{FK}.

 By changing notation and using (\ref{Poh}), one can write (\ref{3.40})
as \be\label{D-det} D|r|^{\a}=b(\a) |r|^{\a-1}, \ee where
\be\label{B1} b(\a)=\a(\a+1/2)\cdots(\a+d-1) \ee is the so-called
{\it Bernstein polynomial} of the determinant \cite{FK}.

\subsection{Bessel functions of matrix argument}
We recall some facts from \cite{Herz} and
 \cite{FK}.

\subsubsection{$\J$-Bessel functions} The $\J$-Bessel function $\J_\nu(r)$,
$r\in\p$, can be defined in terms of the Laplace transform by the
property \be \intl_{\p}
\exp(-\tr(zr))\J_\nu(r)|r|^{\nu-d}dr=\gm(\nu)\exp(-\tr(z^{-1}))\det(z)^{-\nu},
\ee
$$ d=(m+1)/2,\qquad z\in\T_m=\p +i\S_m $$
(the branch of the multi-valued function $\det(z)^{-\nu}$ is chosen so that for $z \in \p$, the argument of $\det(z)^{-\nu}$ is zero).
This gives
\be \J_\nu(r)|r|^{\nu-d}=\frac{\gm(\nu)}{(2\pi i)^N}\intl_{Re\,
z=\sig_0} \exp(\tr(rz- z^{-1}))\det(z)^{-\nu} dz,  \ee
$\sig_0\in\p,  \;  N=m(m+1)/2$, or, by changing variable,

\be \J_\nu(r)=\frac{\gm(\nu)}{(2\pi i)^N}\intl_{Re\, z=\sig_0}
\exp(\tr(z- rz^{-1}))\det(z)^{-\nu} dz. \ee Both integrals are
absolutely convergent for $Re\,\nu>m$. Applying the operator
(\ref{2.50}),  we obtain the following recurrence formulae:
\bea\label{rec1} D
[\J_\nu(r)|r|^{\nu-d}]&=&\frac{\gm(\nu)}{\gm(\nu-1)}\J_{\nu-1}(r)|r|^{\nu-1-d}\;
,\\\label{rec2} D [\J_\nu(r)]&=&(-1)^m
\frac{\gm(\nu)}{\gm(\nu+1)}\J_{\nu+1}(r). \eea The equation
(\ref{rec1}) allows to extend $\J_\nu(r)$ analytically to all
$\nu\in\bbc$.

There is an intimate connection between the Fourier transform and
$\J$-Bessel functions.
\begin{theorem} \cite[p. 492]{Herz}, \cite[p.
355]{FK}
Let $f(x)$ be an integrable  function  on $\Ma$ of the form
$f(x)=f_0(x'x)$, where $f_0$ is a function on $\p$. Then  \be
\label{Boh}(\F f)(y)=\intl_{\Ma} \exp(\tr(iy'x)) f_0 (x'x) dx=
\frac{\pi^{nm/2}}{\gm(n/2)}\tilde f_0 \left(\frac{y'y}{4}\right),
\ee where \be \tilde f_0 (s)= \intl_{\p}\J_{n/2}
(rs)|r|^{n/2-d}f_0(r) dr, \qquad d=(m+1)/2,\ee is the Hankel
transform of $f_0$.
\end{theorem}

This statement is a matrix generalization of the classical result
of S. Bochner, see, e.g., \cite{SW}, Chapter IV, Theorem 3.3.

One can meet another notation for the $\J$-Bessel function in the
literature. We follow the notation from \cite{FK}. It relates to
the notation $A_{\del} (r)$ in \cite{Herz} by the formula
$\J_\nu(r)=\gm(\nu) A_{\del} (r)$, $\del=\nu-d$.  For $m=1$, the
classical Bessel function $J_\nu(r)$ expresses through $\J_\nu(r)$
as
$$
J_\nu(r)=\frac{1}{\Gam(\nu+1)}\left(\frac{r}{2}\right)^\nu
\J_{\nu+1}\left(\frac{r^2}{4}\right).
$$

\subsubsection{$\K$-Bessel functions}\label{sKB} Let, as above,
$r\in\p$, $d=(m+1)/2$. The $\K$-Bessel function $\K_\nu(r)$ is
defined by

\be\label{K1} \K_\nu(r)=\intl_{\p}\exp(\tr(-s-
rs^{-1}))|s|^{\nu-d} ds, \ee or (replace $s$ by $s^{-1}$)

\be\label{K2} \K_\nu(r)=\intl_{\p}\exp(\tr(-s^{-1}-
rs))|s|^{-\nu-d} ds,\ee see \cite{Te}, \cite{FK}.  For $m=1$,
\be \K_\nu(r)=\intl_0^\infty \exp(-s- r/s)s^{\nu-1} ds=2r^{\nu /2}
K_\nu(2\sqrt r), \ee $ K_\nu$ being the classical Macdonald
function.
 \begin{lemma}\label{l1.9}
Let $r\in\p$, $ d=(m+1)/2$.

\noindent{\rm (i)}   The integral (\ref{K1}) (or (\ref{K2}) )
converges absolutely for all $\nu\in\bbc$, and is an entire
function of $\nu$.

\noindent {\rm (ii)}   \ The following estimates hold:

\noindent (a) For $Re \, \nu > d-1$: \be \label{estK1}
 |\K_\nu(r)| \leq  \gm(Re\, \nu);\ee

\noindent (b) For $Re \, \nu < 1-d$: \be \label{estK2}
|\K_\nu(r)|\leq \gm(- Re\, \nu )|r|^{Re\, \nu};\ee

\noindent (c) For $1-d \leq  Re \, \nu \leq d-1$: \be
\label{estK3} |\K_\nu(r)| \leq \gm(d)
+|r|^{1-d-\varepsilon}\gm(d-1+\e), \quad \forall \e>0.\ee

\noindent {\rm (iii)} \ If $Re \, \nu < 1-d$, then
\be\label{K-lim} \lim\limits_{\varepsilon\to 0}
\varepsilon^{-m\nu}\K_\nu(\varepsilon r)=\gm (-\nu)|r|^{\nu} . \ee
\end{lemma}
\begin{proof}
All statements, with  probable exception of (\ref{estK3}), are
known \cite{FK}. The estimate (\ref{estK1}) follows from
(\ref{K1}); (\ref{estK2}) and (\ref{K-lim}) are consequences of
(\ref{K2}). To prove (\ref{estK3}), we write $$ \K_\nu(r)=\left(
\,  \intl_{|s|>1}+ \intl_{|s|<1}\right) \exp(\tr(-s-
rs^{-1}))|s|^{\nu-d} ds=I_1+I_2.$$ For $I_1$ we have
$$
|I_1|<\intl_{|s|>1}\exp(\tr(-s)) ds<\gm(d)
$$
provided $Re \, \nu \leq d$. For $I_2$, by changing variable $s\to
s^{-1}$, we obtain
$$
I_2=\intl_{|s|>1} \exp(\tr(-s^{-1}- rs))|s|^{-\nu-d} ds.
$$
If $Re \, \nu \geq 1-d$, then $|s|^{-Re\, \nu-d}\leq |s|^{\e-1}$
$\forall \e >0$, and \bea\nonumber
 |I_2|&\leq& \intl_{|s|>1}
\exp(\tr(-s^{-1}- rs))|s|^{\e-1} ds\\\nonumber
&<&\intl_{\p}\exp(\tr(-rs))|s|^{\e-1} ds\\\nonumber&=&
\gm(d-1+\e)|r|^{1-d-\varepsilon}. \eea This gives (\ref{estK3}).
\end{proof}

\subsection{Stiefel manifolds}
For $n\geq m$, let $\vnm= \{v \in \frM_{n,m}: v'v=I_m \}$
 be the Stiefel manifold of orthonormal $m$-frames in $\bbr^n$.
 If $n=m$, then $V_{n,n}=O(n)$ is the orthogonal group in $\bbr^n$.
  The group $O(n)$
 acts on $\vnm$ transitively by the rule $g: v\to gv, \quad g\in
 O(n)$, in the sense of matrix multiplication.  The same is true for
  the special orthogonal group  $SO(n)$ provided $n>m$. We fix the corresponding
invariant measure $dv$ on
 $\vnm$
  normalized by \be\label{2.16} \sigma_{n,m}
 \equiv  \! \intl_{\vnm} dv \! = \frac {2^m \pi^{nm/2}} {\gm
 (n/2)}. \ee
\begin{lemma}\label{muir}\cite[p. 589]{Mu} Let $x$ and $y$ be real matrices such that  $x$ is $n \times k$
and $y$ is $m \times k$, $ n \ge m$. Then $x'x=y'y$ if and only if
there exists $v \in V_{n,m}$ such that $x=vy$. In particular, if
$x$ is $n \times m, \; n \ge m$, then there exists $v \in V_{n,m}$
such that $x=vy$, where $y=(x'x)^{1/2}$.
\end{lemma}

\begin{lemma}\label{l2.3} {\rm (polar decomposition).}
Let $x \in \frM_{n,m}, \; n \ge m$. If  $\rank (x)= m$, then \[
x=vr^{1/2}, \qquad v \in \vnm,   \qquad r=x'x \in\p,\] and
$dx=2^{-m} |r|^{(n-m-1)/2} dr dv$.
\end{lemma}

This statement and its generalizations can be found in different
places, see, e.g., \cite[ p. 482]{Herz}, \cite[pp. 66, 591]{Mu},
\cite[p. 130]{FT}. A modification of Lemma \ref{l2.3} in terms of
upper triangular matrices $t\in T_m$ (see (\ref{2.17})) reads as
follows.

\begin{lemma}\label{sph} Let $x \in \frM_{n,m}, \; n \ge m$.
 If  $ \; \rank (x)= m$, then
  \[
x=vt, \qquad v \in \vnm,   \qquad t \in T_m,\] so that
$$
dx=\prod\limits_{j=1}^m t_{j,j}^{n-j} dt_{j,j}dt_{*}dv, \qquad
dt_{*}=\prod\limits_{i<j} dt_{i,j}.
$$
\end{lemma}
\begin{proof}
This statement is also well known. It can be easily derived from
Lemma \ref{l2.3} and (\ref{2.2}). Indeed, if $\rank (x)= m$, then
$x'x\in\p$ and there exists $t\in T_m$ such that $x'x=t't$. We set
$v=xt^{-1}$. Then $v'v=I_m$ and, therefore, $v \in \vnm$. This
proves the representation $x=vt$, where $v \in \vnm$ and  $t \in
T_m$. Furthermore, by Lemma \ref{l2.3} and (\ref{2.2}),
\bea\nonumber \intl_{\Ma} f(x)dx&=&2^{-m}\intl_{\vnm} dv\intl_{\p}
|r|^{(n-m-1)/2}f(vr^{1/2}) dr\\&\stackrel{\rm
(\ref{2.2})}{=}&\intl_{\vnm} dv\intl_{T_m}
f(v(t't)^{1/2})\prod\limits_{j=1}^m t_{j,j}^{n-j}
dt_{j,j}dt_{*}dv\, . \nonumber \eea Now we denote $\lam
=(t't)^{1/2}t^{-1}\in O(m)$, then change the order of integration,
and set $v\lam=u$. This gives
$$
 \intl_{\Ma} f(x)dx=\intl_{T_m}\prod\limits_{j=1}^m t_{j,j}^{n-j} dt_{j,j}dt_{*}
 \intl_{\vnm}f(ut) du,
 $$
 and we are done.
\end{proof}

\section{Radial functions and the Cayley-Laplace operator}

\setcounter{equation}{0}

A function $f(x)$ on $\Ma$  is called radial, if there exists a
function $\f0 (r)$ on $\p$  such that $f(x)=\f0 (x'x)$ for all (or
almost all) matrices $x\in\Ma$. One can readily check that $f$
 is radial if and only if it is $O(n)$ left-invariant, i.e., $f(\gamma
x)=f(x)$ for all $\gamma\in O(n)$.

The {\it Cayley-Laplace operator} $\Del$  on the space $\Ma$ of
matrices $x=(x_{i,j})$ is  defined by \be\label{K-L} \Del=\det(\d
'\d). \ee Here $\partial$ is an $n\times m$  matrix whose entries
are partial derivatives $\d/\d x_{i,j}$. In the Fourier transform
terms, the action of $\Del$ represents a multiplication by the
polynomial $(-1)^m P(y)$, where $y=(y_{i,j})\in\Ma$,
$$
P(y)=|y'y|=\det \left[\begin{array}{cccc} y_1 \cdot y_1 & {.}   & {.}   & y_1 \cdot y_m \\
                              {.} & {.}  & {.}   & {.} \\
                              {.} & {.}   & {.}   & {.} \\
                               y_m  \cdot y_1& {.}     & {.} & y_m \cdot y_m
\end{array} \right],
$$
$ y_1, \dots  y_m$ are column-vectors of the matrix $y$, and
$``\cdot"$ stands for the usual inner product in $\rn$. Clearly,
$P(y)$ is a homogeneous polynomial of degree $2m$ of $nm$
variables $y_{i,j}$, and $\Del$ is a homogeneous differential
operator of order $2m$. For $m=1$, it coincides with the Laplace
operator on $\rn$.

The Cayley-Laplace operator (\ref{K-L}) and its generalizations
were studied by S.P. Khekalo \cite{Kh}. For $m>1$, the operator
$\Del$ is not elliptic because $P(y)=0$ for all non-zero matrices
$y$ of $\rank <m$. Moreover, $\Del$ is not hyperbolic,
 although,
for some $n, m$ and $\ell$, its power $\Del^\ell$ enjoys the
strengthened Huygens' principle; see \cite{Kh} for details.

Our nearest goal is to find  a radial part of $\Del$  corresponding to  the polar decomposition $ x=vr^{1/2}$,
$v \in \vnm$, $r=x'x \in\p$. For $m=1$, the classical result
states  that  the radial part of the Laplace operator on $\rn$ is
$$L=\rho^{1-n}\frac{\d}{\d \rho} \rho^{n-1}\frac{\d}{\d
\rho},\qquad \rho=|x|,\qquad x\in\rn. $$ By changing variable
$r=\rho^2$, we get \be\label{L-rad} L=4r^{1-n/2}\frac{\d}{\d r}
r^{n/2}\frac{\d}{\d r}. \ee The following statement is one of the
main results of the paper. It extends (\ref{L-rad}) to the higher
rank case.

\begin{theorem}\label{TKL-r}
Let $\Omega \subset \Ma$ be an open set consisting of matrices of
rank $m$; $n\geq m\geq 1$. If $f(x)=f_0(x'x)$, $f_0(r)\in
C^{2m}(\p)$, then for  $x \in \Omega$, \be \label{KL-r}(\Del
f)(x)=(Lf_0)(x'x),\ee where \be\label{L-r} L=4^m |r|^{d-n/2}D
|r|^{n/2-d+1}D, \qquad d=(m+1)/2,\ee
$D=\det[((1+\del_{ij})/2)\d/\d r_{i,j}]$ being the operator
(\ref{2.50}).
\end{theorem}
\begin{proof}

For $m=1$, (\ref{L-r}) coincides with (\ref{L-rad}). To prove the
theorem, we first note that without
 loss of generality, one can assume $f$ to be compactly
 supported away from the surface  $\{x: \det(x'x)=0\}$. Otherwise,
 $f(x)$ can be replaced by  $ f_1(x)=\vp(x)\psi(x) f(x)$,
 where $\vp$ and $\psi$ are radial cut-off  functions of the
 form $$\vp(x)=\vp_0(\tr(x'x)), \qquad \psi(x)=\psi_0(\det(x'x)),
 \qquad\vp_0 , \psi_0 \in C^\infty (\bbr_+), $$
$$\begin{array}{ll} \vp_0(\rho) = \left \{
\begin{array} {ll} 1, &
  \mbox{if   $0\leq \rho< N$}, \\
0, & \mbox{if $\rho\geq N+1$},\end{array} \right. \qquad
&\psi_0(\rho) = \left \{
\begin{array} {ll} 0, &
  \mbox{if   $0\leq \rho\leq \e$}, \\
1, & \mbox{if $\rho\geq 2\e$}.\end{array} \right.
\end{array}
 $$
The positive numbers $N$ and $\e$ should be chosen sufficiently
large and small, respectively.

 By the generalized Bochner formula (\ref{Boh}),
$$
\F [\Del f](y)=(-1)^m |y'y|(\F f)(y)=\frac{\pi^{nm/2}}{\gm(n/2)}
h\left(\frac{y'y}{4}\right),
$$
where
$$
h(s)= (-4)^m \intl_{\p}|rs|\J_{n/2} (rs)|r|^{n/2-d-1}f_0(r) dr.
$$
Let us transform $|rs|\J_{n/2} (rs)$.  By (\ref{rec1}) and
(\ref{rec2}), \bea D [|r|^{n/2+1-d}D \J_{n/2}(r)]&=& (-1)^m
\frac{\gm(n/2)}{\gm(n/2+1)} D[\J_{n/2+1}(r)|r|^{n/2+1-d}]\nonumber
\\&=&(-1)^m \J_{n/2}(r)|r|^{n/2-d}. \nonumber \eea Hence,
$$
|r|\J_{n/2} (r)=(-1)^m (\tilde L \J_{n/2}) (r), \qquad \tilde
L=(|r|^{d+1-n/2}D |r|^{n/2-d})(|r|D).
$$
Since $|r|^{d+1-n/2}D |r|^{n/2-d}$  and $|r|D$ are invariant
differential operators with respect to  the transformation $r\to
grg'$, $g\in GL(m,\bbr)$, \cite[p. 294]{FK}, then so is $\tilde
L$. Thus for any $s\in\p$,
$$
|rs|\J_{n/2} (rs)=(-1)^m (\tilde L \J_{n/2}) (sr)=(-1)^m \tilde
L_r[ \J_{n/2} (sr)],
$$
where $\tilde L_r$ stands for the operator $\tilde L $ acting in
the $r$-variable (here we use the symmetry property  $\J_{\nu}
(rs)=\J_{\nu} (s^{1/2}rs^{1/2})$) . It follows that \bea\nonumber
h(s)&=&4^m \intl_{\p}\tilde L_r[ \J_{n/2} (sr)]|r|^{n/2-d-1}f_0(r)
dr\\\nonumber &=&4^m \intl_{\p}(D_{r}  |r|^{n/2-d+1}D_{r} )[
\J_{n/2} (sr)]f_0(r) dr .\eea
 Owing to remark at the beginning of the proof, one can integrate by
 parts and get
 \bea\nonumber h(s)&=&(-4)^m \intl_{\p}D_{r}[ \J_{n/2} (sr)] |r|^{n/2-d+1}(D f_0)(r)
dr\\\nonumber &=&\intl_{\p} \J_{n/2} (sr) |r|^{n/2-d}(Lf_0)(r) dr
 ,\eea
 where
 $L=4^m |r|^{d-n/2}D |r|^{n/2-d+1}D$. Thus
 $$
 \F [\Del f](y)=\frac{\pi^{nm/2}}{\gm(n/2)}\intl_{\p} \J_{n/2} \left(\frac{1}{4}r y'y\right) |r|^{n/2-d}(Lf_0)(r)
 dr,
 $$
 which implies (\ref{KL-r}) .
\end{proof}

\begin{example}
Let $f(x)=|x|_m^\lam, \; |x|_m=\det (x'x)^{1/2}$. By Theorem
\ref{TKL-r} and (\ref{D-det}), $(\Del f)(x)=\vp(x'x)$, where \bea
\nonumber\vp(r)&=&4^m |r|^{d-n/2}D
|r|^{n/2-d+1}D|r|^{\lam/2}\\[14pt]\nonumber&=& 4^m b(\lam/2)
|r|^{d-n/2}D |r|^{(n+\lam)/2-d}\\[14pt]\nonumber&=& 4^m
b(\lam/2) b((n+\lam)/2-d) |r|^{\lam/2-1}.
 \eea
Thus, we have arrived at the following identity of the Bernstein
type \be \label{Dxm}\Del |x|_m^\lam =\B(\lam) |x|_m^{\lam -2},\ee
where, owing to (\ref{B1}), the polynomial $\B(\lam)$ has the form
\be\label{B} \B(\lam)=(-1)^m
\prod\limits_{i=0}^{m-1}(\lam+i)(2-n-\lam+i). \ee

An obvious consequence of (\ref{Dxm}) in a slightly different
notation reads \be\label{vaz}\Del ^k
|x|_m^{\a+2k-n}=B_k(\a)|x|_m^{\a-n},\ee  \bea\label{bka}
B_k(\a)&=&\prod\limits_{i=0}^{m-1}\prod\limits_{j=0}^{k-1}(\a-i+2j)(\a-n+2+2j+i)
\\ &=& B_k(n-\a-2k). \nonumber \eea
\end{example}

\section{Zeta integrals}

\subsection{Definition and example} Let us consider the  zeta   integral \be\label{zeta}
\Z(f,\a-n)=\intl_{\Ma} f(x) |x|^{\a-n}_m dx  \ee where $f(x)$ is a
Schwartz function on
 $\Ma$,
$n\geq m$, $|x|_m =\det (x'x)^{1/2}$. The following example gives
a flavor of basic properties of $\Z(f, \a-n)$.
\begin{example}
 Let $e(x)=\exp(-\tr(x'x))
$ be  the Gaussian function. By Lemma \ref{l2.3}, for $Re\, \a
>m-1$ we have \bea \Z(e,
\a-n)&=&2^{-m}\sig_{n,m}\intl_{\p}|r|^{\a/2 -d}
\exp(-\tr(r))dr\\\nonumber &=&c_{n,m}\gm(\a/2),\quad
c_{n,m}=\frac{\pi^{nm/2}}{\gm(n/2)}. \eea $d=(m+1)/2$. On the
other hand, the well known formula for the Fourier transform
yields \be (\F e)(y)= \pi^{nm/2}e(y/2),\ee and for $Re\, \a
<n-m+1$ we obtain \bea \nonumber \Z(\F e, -\a)&=&c_{n,m}
\pi^{nm/2}\intl_{\p}|r|^{(n-\a)/2 -d} \exp(-\tr(r/4))dr\\
&=& d_{n,m} \gm((n-\a)/2), \quad d_{n,m}= c_{n,m} \pi^{nm/2}
2^{m(n-\a)}. \eea Thus, after analytic continuation we obtain the
following meromorphic functions: \be \label{ac-ze1}  \Z(e, \a-n)=
c_{n,m}\gm(\a/2),\qquad \quad \;\; \a\neq m-1,m-2, \dots\; ,\ee
\be\label{ac-ze2} \Z(\F e, -\a)= d_{n,m} \gm((n-\a)/2),\quad
\a\neq n-m+1,n-m+2,\dots  \, .
 \ee
These equalities imply the functional relation \be\label{rel-z}
\frac{\Z(e, \a-n)}{\gm(\a/2)}=\pi^{-nm/2} 2^{m(\a-n)}\frac{\Z(\F
e, -\a)}{\gm((n-\a)/2)}\; , \ee which is a prototype of similar
formulas for much more general zeta functions.  Since each side
 of  (\ref{rel-z}) equals  $c_{n,m}$, this formula extends to
all complex $\a$. Note that excluded values of $\a$ in
(\ref{ac-ze1}) and (\ref{ac-ze2}) correspond to $m\geq 2$. If
$m=1$, they proceed with step $2$, namely, $\a\neq 0,-2,-4, \dots
,\; $ and $\a\neq n,n+2, \dots\; $, respectively.
\end{example}

In the following, throughout the paper, we assume $m\geq 2$.

\begin{lemma}\label{lacz}
Let $f\in\S(\Ma)$. For $Re\, \a > m-1$,  the integral (\ref{zeta})
is absolutely convergent, and for $Re\, \a \leq m-1$, it extends
as a meromorphic function of $ \; \a$ with the only poles $ \;
m-1, m-2,\dots\;$. These poles and their orders are exactly the
same as of the gamma function $\gm(\a/2)$. The normalized  zeta
integral $ \Z(f,\a-n)/\gm(\a/2)$
 is an entire
function of $\a$.
\end{lemma}
\begin{proof}
This statement is known; see, e.g., \cite{Kh}, [Sh]. We present
the proof for the sake of completeness. The equality
(\ref{ac-ze1}) says that the function $\a \to \Z(f,\a-n)$ has
poles at least at the same points and  at least of the same order
as the gamma function $\gm(\a/2)$. Our aim is to show that no
other poles occur, and the orders cannot exceed  those of
 $\gm(\a/2)$.
Let us transform (\ref{zeta}) by passing to upper triangular
matrices $t\in T_m$ according to Lemma \ref{sph}. We have
\be\label{z-tr} \Z(f,\a-n)=\intl_{\bbr^m_+} F(t_{1,1},\dots ,
t_{m,m}) \prod\limits_{i=1}^m t_{i,i}^{\a-i} dt_{i,i}\; , \ee
$$
F(t_{1,1},\dots , t_{m,m}) =\!\!\!\!\!\! \intl_{\bbr^{m(m-1)/2}}
\!\!\!\! dt_{*} \intl_{\vnm} f(vt) dv, \quad
dt_{*}=\prod\limits_{i<j} d t_{i,j}.
$$
Since $F$ extends as an even Schwartz function in each argument,
it can be written as
$$
F(t_{1,1},\dots , t_{m,m}) =F_0(t^2_{1,1},\dots , t^2_{m,m}),
$$
where $F_0\in\S(\bbr^m)$ (use, e.g., Lemma 5.4 from \cite[p.
56]{Tr}. Replacing $t_{i,i}^2$ by $s_{i,i}$, we represent
(\ref{z-tr}) as a direct product of one-dimensional distributions
\be\label{z-h} \Z(f,\a-n)=2^{-m}(\prod\limits_{i=1}^m
(s_{i,i})_+^{(\a-i-1)/2},\;  F_0(s_{1,1},\dots , s_{m,m})),\ee
which is a meromorphic function of $\a$ with the poles $m-1, m-2,
\dots \;$, see \cite{GSh1}. These poles and their orders coincide
with those  of the gamma function $\gm(\a/2)$, cf. (\ref{ac-ze1}).
To normalize the function (\ref{z-h}), following \cite{GSh1}, we
divide it by the product
$$\prod\limits_{i=1}^m\Gam((\a-i+1)/2)=\prod\limits_{i=0}^{m-1}\Gam((\a-i)/2)=
\gm(\a/2)/\pi^{m(m-1)/4}.$$ As a result we obtain an entire
function.
\end{proof}

\subsection{A functional equation for the zeta integral}
Let us prove another main result of the paper.
\begin{theorem}\label{tZ-F}
 If $f\in\S(\Ma)$, $n\geq m$, then
 \be\label{z-f}
\frac{\Z(f, \a-n)}{\gm(\a/2)}=\pi^{-nm/2} 2^{m(\a - n)}\frac{\Z(\F
f, -\a)}{\gm((n-\a)/2)}. \ee
\end{theorem}
\begin{proof}
We recall  that both sides of this equality  are understood in the
sense of analytic continuation. Furthermore, the distribution on
the left-hand side of (\ref{z-f}) is regular if and only if $Re\,
\a >m-1$, whereas the right-hand side is regular if and only if
$Re\, \a <n-m+1$. For $2m\geq n+2$ these two sets have no common points.
 We consider the cases $2m<n+2$ and $2m\geq n+2$ separately.

{\rm (i)} {\bf The case  $2m<n+2$.} We start with the following
equality: \be\label{3.5}
\begin{array} {ll}\displaystyle{\intl_{\Ma} (\F f)(y)|\e
I_m+y'y|^{-\a/2}dy}\\[14pt]
=\displaystyle{\frac{\pi^{nm/2}\e^{m(n-\a)/2}}{\gm(\a/2)}
\intl_{\Ma} f(x) \K_{(\a-n)/2}\left(\frac{\e}{4}x'x\right)
dx},\qquad \e>0.
 \end{array} \ee
 Here $\K_{(\a-n)/2}$ is the $\K$-Bessel function (see Section
\ref{sKB}),  and \be\label{d1} m-1< Re\, \a <n-m+1\ee (since
$2m<n+2$, the domain (\ref{d1}) is not void). Suppose for a
moment, that (\ref{3.5}) is true. Then the result follows if we
pass to the limit as $\e\to 0$. Indeed, by (\ref{estK2})  for all
$\e>0$ we have \be\label{est-Ke}
|\e^{m(n-\a)/2}\K_{(\a-n)/2}\left(\frac{\e}{4}x'x\right) |\leq
\gm\left(\frac{n-Re\, \a}{2}\right)\Big |\frac{x'x}{4}\Big
|^{(Re\, \a-n)/2}. \ee Furthermore, by (\ref{K-lim}),
\be\label{lim-Ke} \lim\limits_{\e\to
0}\e^{m(n-\a)/2}\K_{(\a-n)/2}\left(\frac{\e}{4}x'x\right)=2^{m(n-\a)}\gm\left(\frac{n-
\a}{2}\right)|x|_m^{\a-n}.\ee Hence, by the Lebesgue theorem on
dominated convergence, we are done.

To prove (\ref{3.5}), let $e_s (x)=\exp (-\tr(xsx')/4\pi), \; s
\in \p$. By the Plancherel formula, \be\label{pl1}
|s|^{-n/2}\intl_{\Ma} (\F f)(y) \exp (-\tr (\pi y
s^{-1}y'))dy=\intl_{\Ma}  f(x) e_s (x) dx. \ee  Now we multiply
(\ref{pl1}) by $|s|^{(n-\a)/2 -d}\exp (-\tr(\e\pi s^{-1})), \; d=
(m+1)/2$, then integrate in $s$, and change the order of
integration. We obtain
\[ \intl_{\Ma} (\F f)(y) a (y) dy=\intl_{\Ma}  f(x) b(x)
dx, \] where \bea a (y)&=&\intl_{\p} |s|^{-\a/2 -d} \exp [-\tr
(\pi s^{-1}(y'y+\e I_m))] ds \qquad (s=t^{-1})\nonumber
\\\nonumber &=& \intl_{\p} |t|^{\a/2 -d} \exp [-\tr (\pi t(y'y+\e
I_m))] dt\\&=&\Gam_m (\a/2)\pi ^{-\a m/2}|y'y+\e I_m|^{-\a/2}
,\qquad Re \,\a > m-1, \nonumber \eea  and \bea \nonumber
b(x)&=&\intl_{\p} |s|^{(n-\a)/2 -d} \exp [-\tr ( sx'x/4\pi+\e \pi
s^{-1} ))]  ds\\\nonumber &=&(\e \pi)^{m(n-\a)/2} \intl_{\p}
|u|^{(n-\a)/2 -d} \exp [-\tr ( u \e x'x/4+u^{-1} ))] du\\
\nonumber &=& (\e \pi)^{m(n-\a)/2} \, \K_
{(\a-n)/2}\left(\frac{\e}{4}x'x\right). \eea This gives
(\ref{3.5}). Application of the Fubini theorem  in this argument
is justified because both integrals in (\ref{3.5}) converge
absolutely. Indeed, for the left-hand side the absolute
convergence is obvious, and for the right-hand side it follows
from (\ref{lim-Ke}).

{\rm (ii)} {\bf The case  $2m\geq n+2$.} In this case the interval
(\ref{d1}) is void. To circumvent this difficulty, we replace
$f(x)$ in (\ref{pl1}) by $|x|_m^{2k}\Del ^k f(x)$ and proceed as
above, assuming $k$ large enough and $Re \,\a > m-1$. Formally we
obtain \be\label{pl2}
\begin{array}{ll}\displaystyle{\intl_{\Ma}\Del ^k [|y|_m^{2k}(\F
f)(y)] \, |\e
I_m+y'y|^{-\a/2}dy}\\=\displaystyle{\frac{\pi^{nm/2}\e
^{m(n-\a)/2} }{\Gam_m (\a/2)}\intl_{\Ma}  |x|_m^{2k}\Del ^k
f(x)\K_{(\a-n)/2}\left(\frac{\e}{4}x'x\right)  dx}. \end{array}\ee
The equality (\ref{pl2}) will become meaningful if we justify
application of the Fubini theorem and specify a suitable interval
for $Re\,\a$. Clearly, the left-hand side of (\ref{pl2})
absolutely converges {\it for all $\a\in\bbc$}, and one should
take care of the right-hand side only. By Lemma \ref{l1.9}
(\rm{ii}), the integral on the right-hand side absolutely
converges for $Re \,\a
> m-1-2k$ provided $2k>2m-2-n$. Indeed, for small $|x|_m$ we have
$$
|x|_m^{2k}\Del ^k
f(x)\K_{(\a-n)/2}\left(\frac{\e}{4}x'x\right)=O(|x|_m^{\lam -n}),
$$
where $$ \lam= \left \{
\begin{array} {ll} 2k+n, &
 \mbox{ if   $Re \,\a >n+ m-1$}, \\
  2k+n+1-m-\e_0,
\forall \e_0>0, & \mbox { if   $n-m+1\leq Re \,\a \leq n+ m+1$}, \\
2k+Re \,\a, & \mbox { if   $Re \,\a <n- m+1$}.
\end{array}
\right.
 $$
 If $Re \,\a > m-1-2k$ and $\e_0$ is
small enough ($\e_0<2k-2m+n+2$) then  $\lam
>m-1$  and  the desired convergence
follows. Thus, (\ref{pl2}) is justified for $Re \,\a > m-1$. Since
both sides of (\ref{pl2}) are analytic functions of $\a$ in a
larger domain $Re \,\a > m-1-2k$ containing a non-void strip
\be\label{d2}  m-1-2k<Re \,\a < n-m+1,\ee   (\ref{pl2}) also holds
in this strip, and one can utilize (\ref{est-Ke}) and
(\ref{lim-Ke}) in order to pass to the limit as $\e\to 0$. This
gives \be\label{3.6}  \Z(\Del ^k [ |y|_m^{2k} (\F f)(y)],
-\a)=c_\a \Z(\Del ^k f, \a+2k-n), \ee  \be\label{ca} c_\a=
\pi^{nm/2}
2^{m(n-\a)}\gm\Big(\frac{n-\a}{2}\Big)/\gm\Big(\frac{\a}{2}\Big).\ee
To transform (\ref{3.6}), we use the equality  \be\Del ^k
|x|_m^{\a+2k-n}=B_k(\a)|x|_m^{\a-n},\ee  $$
B_k(\a)=\prod\limits_{i=0}^{m-1}\prod\limits_{j=0}^{k-1}(\a-i+2j)(\a-n+2+2j+i)=
B_k(n-\a-2k),$$ see (\ref{vaz}). Then
  the right-hand side of (\ref{3.6}) becomes $c_\a
  B_k(\a)\Z(f,\a-n)$. For the left-hand side we obtain (set $\vp(y)=|y|_m^{2k} (\F
f)(y)$, $ \; \a=n-\b-2k)$ \bea \Z(\Del ^k[
|y|_m^{2k} (\F f)(y)], -\a)&=& \Z(\Del ^k \vp, \b+2k-n) \nonumber \\
&=&B_k(\b)\Z( \vp, \b-n) \nonumber \\&=&B_k(n-\a-2k)\Z( \F f,
-\a)\nonumber \\
&=&B_k(\a)\Z( \F f, -\a). \nonumber \eea
 Finally, (\ref{3.6}) reads
$$
\Z( \F f, -\a)=c_\a \Z(  f, \a-n), \qquad  m-1-2k<Re \,\a < n-m+1,
$$
and the desired equality (\ref{z-f}) follows by analytic continuation.
\end{proof}

\subsection{Normalized zeta distributions of integral order}

It is convenient to introduce a special notation for the
normalized zeta integral $\Z(f,\a-n)/\gm(\a/2)$ which is
  an entire function of $\a$. We denote \be\label{nzd} \z_\a (x) \!
=\!\frac{|x|_m^{\a -n}}{\gm(\a/2)}, \quad (\z_\a, f)=a.c.
\frac{1}{\gm(\a/2)}\intl_{\frM_{n,m}} \! f(x)|x|_m^{\a -n} dx,\ee
where``$a.c.$'' abbreviates analytic continuation. We call $\z_\a$
{\it a normalized zeta distribution of order $\a$}.

In view of forthcoming applications in Section
6, normalized zeta distributions of integral
order deserve special treatment.
 A striking feature  of the
distribution $\z_\a$ is that for $\a=k, \; k=0,1,2, \dots m-1 $,
(outside of the domain of absolute convergence !) it is a positive
measure. This measure is supported by a lower-dimensional manifold
(in the rank-one case $m=1$ we have only one point $\a=0$
corresponding to the delta function at the origin).
\begin{definition} The set
\be\label{wal} \W =\{0,1,2, \ldots,  m-1\} \cup \{\a: Re \, \a >
m-1 \}\ee will be called {\it a Wallah set} of the normalized zeta
distribution $\z_\a$.
\end{definition}

This notion was introduced in \cite{FK} for Riesz distributions
associated to symmetric cones. Below we  obtain explicit
representations of $\z_\a$ for integral values of $\a$
(specifically, for $0<\a \le n$) in the Wallah set (\ref{wal}).
\begin{theorem}\label{tzk}  Let $f  \in  \S(\Ma)$. For
 $\a=k$, $k=1, 2,
\ldots,  n$, \be\label{znk}(\z_k,
f)=\frac{\pi^{(n-k)m/2}}{\gm(n/2)} \intl_{SO(n)} d\gam
\intl_{\frM_{k,m}}f \left (\gam \left[\begin{array} {c} \om \\ 0
\end{array} \right]  \right ) \, d\om. \ee
Furthermore, in the case
 $\a=0$ we have \be\label{zn0} (\z_0, f)=\frac{\pi^{nm/2}}{\gm(n/2)}f(0). \ee
\end{theorem}
\begin{proof} STEP 1.
Let first $k>m-1$ . In polar coordinates we have
 \bea\nonumber
\Z(f,k-n)&=&\intl_{\Ma} f(x) |x|^{k-n}_m dx\\\nonumber &=&
2^{-m}\sig_{n,m}\intl_{\P_{m}}|r|^{k/2 -d} dr
\intl_{SO(n)} f\left ( \gam \left[\begin{array} {c} r^{1/2} \\
0\end{array} \right]\right ) d\gam .\eea Now we replace $\gam$ by
$\gam \left[\begin{array} {cc} \b& 0 \\
0& I_{n-k}\end{array} \right]$, $\b\in SO(k)$, then integrate in
$\b\in SO(k)$, and replace the integration over $SO(k)$ by that
over $V_{k,m}$. We get

\bea\nonumber \Z(f,k-n)&=&
\frac{2^{-m}\sig_{n,m}}{\sig_{k,m}}\intl_{SO(n)} d\gam
\intl_{\P_{m}}|r|^{k/2 -d} dr \intl_{V_{k,m}}
 f\left ( \gam \left[\begin{array} {c} v r^{1/2} \\
0\end{array} \right]\right ) dv \\\nonumber
&{}&({\mbox set \quad  \om=v r^{1/2}\in \Mkm})\\
\nonumber &=& \frac{\sig_{n,m}}{\sig_{k,m}}\intl_{\frM_{k,m}}d\om
\intl_{SO(n)}f \left (\gam \left[\begin{array} {c} \om \\ 0
\end{array} \right]  \right ) \, d\gam.\eea
This coincides with (\ref{znk}).

 STEP 2.
Our next task is to prove that analytic continuation of $(\z_\a,
f)$ at the point $\a=k$ ($\leq m-1$) has the form (\ref{znk}). We
first note that for $\a=0$, (\ref{zn0}) is an immediate
consequence of (\ref{z-f}).  Let $k>0$.
 We split $x \in
\Ma$ in two blocks $x=[y; b]$ where $ y \in \frM_{n,k}$ and $ b
\in \frM_{n, m-k}$. Then for $Re \, \a>m-1,$
\[ (\z_\a, f)=\frac{1}{\gm(\a/2)} \intl_{\frM_{n,k}} dy
\intl_{\frM_{n,m-k}}f([y; b])\left |\begin{array}{ll}  y'y & y'b \\
b'y& b'b
\end{array}\right|^{(\a-n)/2} db \]
where $\left |\begin{array}{ll}  * & * \\
*&*
\end{array}\right|$ denotes the determinant of the respective
matrix $\left [\begin{array}{ll}  * & * \\
*&*
\end{array}\right]$. By passing to polar coordinates (see Lemma \ref{l2.3})
$y=vr^{1/2}, \; v \in V_{n,k}$, $ r \in \P_k$, we have \bea
(\z_\a, f)&=&\frac{2^{-k}}{\gm(\a/2)} \intl_{V_{n,k}} dv
\intl_{\P_k} |r|^{(n-k-1)/2} dr \nonumber \\ &\times&
\intl_{\frM_{n,m-k}}f([vr^{1/2}; b])\left |\begin{array}{ll}  r & r^{1/2}v'b \\
b'vr^{1/2}& b'b
\end{array}\right|^{(\a-n)/2} db \nonumber \\
&=&\frac{2^{-k} \,\sig_{n,k} }{\gm(\a/2)} \intl_{SO(n)}d\gam
\intl_{\P_k} |r|^{(n-k-1)/2} dr  \nonumber \\ &\times&
\intl_{\frM_{n,m-k}}f_\gam ([\lam_0 r^{1/2}; b])\left |\begin{array}{ll}  r & r^{1/2}\lam'_0 b \\
b'\lam_0 r^{1/2}& b'b
\end{array}\right|^{(\a-n)/2} db. \nonumber \eea
Here \[ \lam_0= \left[\begin{array} {c} I_k \\ 0
\end{array} \right]  \in V_{n,k}, \qquad f_\gam (x)=f (\gam x).\]
We write \[ b=\left[\begin{array} {c} b_1 \\ b_2
\end{array} \right], \qquad b_1 \in \frM_{k,m-k}, \qquad b_2 \in
\frM_{n-k,m-k}.\] Since $ \lam'_0 b=b_1$, then \bea
(\z_\a,f)&=&\frac{2^{-k} \,\sig_{n,k} }{\gm(\a/2)}
\intl_{SO(n)}d\gam \intl_{\P_k} |r|^{(n-k-1)/2} dr
\intl_{\frM_{k,m-k}} db_1 \nonumber \\
&\times& \intl_{\frM_{n-k,m-k}}f_\gam \left ( \left [\begin{array}{ll}   r^{1/2}&  b_1 \\
0 & b_2
\end{array}\right]\right )\, \left |\begin{array}{ll}  r & r^{1/2} b_1 \\
b'_1 r^{1/2}& b'_1 b_1 +b'_2 b_2
\end{array}\right|^{(\a-n)/2} db_2. \nonumber \eea
Note that \[ \left [\begin{array}{ll}  r & r^{1/2} b_1 \\
b'_1 r^{1/2}& b'_1 b_1 +b'_2 b_2
\end{array}\right]=\left [\begin{array}{ll}  r & 0 \\
b'_1 r^{1/2}& I_{m-k}
\end{array}\right] \left [\begin{array}{ll}  I_{k}& r^{-1/2} b_1 \\
0 & b'_2 b_2
\end{array}\right], \] and  \[ \det\left [\begin{array}{ll}  r & r^{1/2} b_1 \\
b'_1 r^{1/2}& b'_1 b_1 +b'_2 b_2
\end{array}\right ]=\det (r) \det (b'_2 b_2),\] see, e.g., \cite[ p.
577]{Mu}. Therefore, \be \label{an}(\z_\a,f)=c_\a
\intl_{SO(n)}d\gam \intl_{\P_k} |r|^{(\a-k-1)/2} dr
\intl_{\frM_{k,m-k}} \psi_{\a -k} (\gam, r, b_1) db_1,\ee where
\[ c_\a=\frac{2^{-k} \, \sig_{n,k} \, \Gam_{m-k} ((\a -k)/2)}{\gm
(\a/2)}, \] \bea  \psi_{\a -k} (\gam, r,
b_1)&=&\frac{1}{\Gam_{m-k} ((\a -k)/2)}
\intl_{\frM_{n-k,m-k}}f_\gam \left ( \left [\begin{array}{ll}   r^{1/2}&  b_1 \\
0 & b_2
\end{array}\right]\right ) \nonumber \\ &\times&  |b'_2 b_2|^{(\a-k)/2-(n-k)/2} db_2. \nonumber
\eea The last expression represents the normalized zeta
distribution of order $\a-k$ in the $b_2$-variable. Owing to
(\ref{zn0}), analytic continuation of (\ref {an}) at $\a=k$ reads
\[(\z_k,f)=c_k \intl_{SO(n)}d\gam
\intl_{\P_k} |r|^{-1/2} dr \intl_{\frM_{k,m-k}} f_\gam \left ( \left [\begin{array}{ll}   r^{1/2}&  b_1 \\
0 & 0 \end{array}\right]\right ) db_1,\]
$c_k=[\pi^{(n-k)(m-k)/2}/\Gam_{m-k}((n-k)/2)]\lim\limits_{\a \to
k}c_\a$. To transform this expression, we replace $\gam$ by $\gam \left [\begin{array}{ll}   \b& 0 \\
0 & I_{n-k} \end{array}\right]$, $ \b \in SO(k)$, and integrate in
$\b$. This gives \bea (\z_k,f)&=&c_k   \intl_{\P_k} |r|^{-1/2} dr
\intl_{SO(k)}  d\b   \intl_{\frM_{k,m-k}}
db_1 \nonumber \\
&\times&\intl_{SO(n)}   f_\gam \left ( \left [\begin{array}{ll}   \b r^{1/2}&  \b b_1 \\
0 & 0
\end{array}\right]\right ) d\gam \nonumber \\
&& \text{\rm (set $\zeta=\b b_1, \quad \eta=b |r|^{1/2}$ and use
Lemma \ref{l2.3})}
\nonumber \\
&=&\frac{2^{k} \, c_k}{\sig_{k,k}}\intl_{\frM_{k,k}} d\eta
\intl_{\frM_{k,m-k}} d\zeta\intl_{SO(n)} f_\gam \left ( \left [\begin{array}{ll}   \eta &  \zeta \\
0 & 0
\end{array}\right]\right ) d\gam \nonumber \\
&=&c\intl_{SO(n)} d\gam \intl_{\frM_{k,m}}f \left (\gam
\left[\begin{array} {c} \om \\ 0
\end{array} \right]  \right ) \, d\om, \nonumber \eea
\[ c=\frac{\pi^{(n-k)(m-k)/2}\sig_{n,k}}{\sig_{k,k}\Gam_{m-k}((n-k)/2)}  \lim\limits_{\a \to
k}  \frac{ \Gam_{m-k} ((\a -k)/2)}{\gm
(\a/2)}=\frac{\pi^{(n-k)m/2}}{\gm(n/2)}.\] (here we used formulae
(\ref{2.5.1}) and (\ref{2.16})). \end{proof}

The following formulas are consequences of  (\ref{znk}).
\begin{corollary}\label{nf}
For all  $k=1, 2, \ldots,  n$, \be\label{r1}
(\z_k,f)=c_1\intl_{V_{n,k}} dv\intl_{\frM_{k,m}} f(v \om)d\om,\ee
\be\label{con1} c_1= 2^{-k} \, \pi^{(nm-km-nk)/2} \,
\Gam_k(n/2)/\Gam_m (n/2). \ee Moreover, if $k=1, 2, \ldots, m-1$,
then \be\label{r3}(\z_k,f)=c_1\intl_{V_{m,k}} du
\intl_{\frM_{n,k}}f(yu') |y|_k^{m-n} dy, \ee and \be\label{r2}
(\z_k,f)=c_2\intl_{\frM_{n,k}}
\frac{dy}{|y|_k^{n-m}}\intl_{\frM_{k,m-k}}  f([y; yz]) dz, \ee
\be\label{con2} c_2= \pi^{(m-k)(n/2-k)}\,/\Gam_k (k/2) \,
\Gam_{m-k}((n-k)/2). \ee
\end{corollary}
\begin{proof} From (\ref{znk}) we have
\bea  (\z_k, f)&=&\frac{\pi^{(n-k)m/2}}{\gm(n/2)}
\intl_{\frM_{k,m}}d\om \intl_{SO(n)}f(\gam\lam_0 \om)d\gam \quad
\left ( \lam_0= \left[\begin{array} {c} I_k \\ 0
\end{array} \right]  \in V_{n,k} \right )\nonumber \\
&=&\frac{\pi^{(n-k)m/2}}{\sig_{n,k} \,
\gm(n/2)}\intl_{\frM_{k,m}}d\om \intl_{V_{n,k}}f(v \om)dv\nonumber
\eea which coincides with (\ref{r1}).  To prove (\ref{r3}), we
pass to polar coordinates  in (\ref{r1}) by setting
$\om'=ur^{1/2}, \; u \in V_{m,k}, \; r \in \P_k$. This gives \bea
(\z_k, f)&=& 2^{-k} \, c_1
\intl_{V_{n,k}}dv\intl_{\P_k}|r|^{(m-k-1)/2}
dr\intl_{V_{m,k}} f(vr^{1/2}u') du\nonumber \\
&=&c_1\intl_{V_{m,k}} du \intl_{\frM_{n,k}}f(yu') |y|_k^{m-n}
dy.\nonumber \eea To prove (\ref{r2}), we represent $\om$ in
(\ref{r1}) in the block form $u=[\eta; \zeta], \; \eta \in
\frM_{k,k}, \; \zeta \in \frM_{k,m-k}$, and change the variable
$\zeta=\eta z$. This gives
\[ (\z_k, f)=c_1\intl_{\frM_{k,k}}
|\eta|^{m-k}d\eta\intl_{\frM_{k,m-k}}dz \intl_{V_{n,k}}f(v[\eta;
\eta z])dv.\] Using Lemma \ref{l2.3} repeatedly, and changing
variables, we obtain \bea  (\z_k, f)&=&2^{-k} \, c_1  \,
\sig_{k,k}\intl_{\P_k} |r|^{(m-k-1)/2}
dr\intl_{\frM_{k,m-k}}dz \intl_{V_{n,k}}f(v[r^{1/2}, r^{1/2}z])dv\nonumber \\
&=&c_1 \, \sig_{k,k}\intl_{\frM_{n,k}}
\frac{dy}{|y|_k^{n-m}}\intl_{\frM_{k,m-k}} f([y; yz]) dz.\nonumber
\eea By  (\ref{con1}), (\ref{2.16}) and (\ref{2.5.2}), this
coincides with (\ref{r2}). \end{proof}

The representation (\ref{r2})
 was obtained in \cite{Sh} and \cite{Kh} in a different way.
 An idea of (\ref{r3}) is due to E. Ournucheva.
\begin{remark}
One can also write $ (\z_k,f)$  as \be\label{z-mes}
(\z_k,f)=\intl_{\Ma} f(x)d\nu_k(x),\qquad f\in\S(\Ma), \ee where
$\nu_k$ is a positive locally finite measure defined by
\be\label{mes} (\nu_k ,\vp)\equiv c_1\intl_{V_{n,k}}
dv\intl_{\frM_{k,m}} \vp(v \om)d\om, \quad \vp \in C_c(\Ma),
 \ee
$C_c(\Ma)$ being the space of compactly supported continuous
functions on $\Ma$; cf. (\ref{r1}).
 In order to characterize the support of $\nu_k$, we denote \be
\frM_{n,m}(k)  =  \{x: x   \in   \Ma, \, \rank (x)  =  k \}, \ee
\be \bar \frM_{n,m}(k)  = \bigcup_{j=0}^{k} \frM_{n,m}(j) \qquad \text{\rm 
(the closure of $ \frM_{n,m}(k)$)}. \ee
\begin{lemma}\label{nu}
The following statements hold.\\
 \noindent{\rm(i)}  $\supp\,
\nu_k=\bar \frM_{n,m}(k)$.

\noindent{\rm(ii)} The manifold $ \; \frM_{n,m}(k)$ is an orbit of
 $\;  e_k=\left[\begin{array}{ll}  I_k&0\\
0&0
\end{array}\right]_{n \times m}$
 under the group of transformations $$x \to g_1xg_2, \qquad  g_1 \in
 GL(n,\bbr), \quad   g_2 \in
 GL(m,\bbr).$$

\noindent{\rm(iii)} The manifold $ \; \bar \frM_{n,m}(k)$ is a
collection of all matrices $x \in \Ma$ of the form
\be\label{form1} x=\gam \left[\begin{array} {c} \om \\ 0
\end{array} \right], \quad \g \in SO(n), \quad \om \in
\frM_{k,m}, \ee or \be\label{form2} x=v\om, \qquad v \in V_{n,k},
\quad \om \in \frM_{k,m}. \ee
\end{lemma}
\begin{proof} (i) Let us consider (\ref{mes}).
Since $\rank (v \om)\leq k$, then $(\nu_k ,\vp)=0$ for all $\vp
\in C_c(\Ma)$  supported away from $\bar \frM_{n,m}(k)$. This
means that $\supp\, \nu_k=\bar \frM_{n,m}(k)$.

(ii) Let us show that each $x\in \frM_{n,m}(k)$ is represented in
the form $x=g_1 e_k g_2$ for some $g_1 \in
 GL(n,\bbr)$ and $ g_2 \in
 GL(m,\bbr)$. By Lemma \ref
{muir}, each matrix $x$ can be written as $x=u\rho$, where
$u\in\vnm$ and $\rho=(x'x)^{1/2}$ is a positive semi-definite
$m\times m$ matrix of rank $k$. By taking into account that
$\rho=g'_2 \left[\begin{array}{ll}  I_k&0\\
0&0
\end{array}\right]_{m\times m} g_2 \, $ for some $ g_2 \in
 GL(m,\bbr)$, and $u=\gam
\left[\begin{array} {c} I_m
\\ 0
\end{array} \right]$
for some $\g\in O(n)$, we obtain $x=g_1 e_k g_2$ with
$$ g_1=
\g\left[\begin{array}{ll}  g_2'&0\\
0& I_{n-m}
\end{array}\right]\in
 GL(n,\bbr).$$

(iii) Clearly, each matrix of the form (\ref{form1}) or
(\ref{form2}) has rank $\le k$. Conversely, if $\rank (x) \le k$
then, as above, $x=u\rho =\gam \left[\begin{array} {c} \rho \\ 0
\end{array} \right]$ where $ \g \in O(n)$ and $\rho=(x'x)^{1/2}$ is a positive semi-definite
$m\times m$ matrix of rank $\le k$. The latter can be written as
\[ \rho=g\lam g', \qquad  g \in O(m), \qquad \lam =\diag(\lam_1,
\ldots , \lam_k, 0, \ldots , 0),\] and therefore, \[ x=\gam
\left[\begin{array} {c} g\lam \\ 0
\end{array} \right]g'=\g\left[\begin{array}{ll}  g&0\\
0& I_{n-m}
\end{array}\right]\left[\begin{array} {c} \lam \\ 0
\end{array} \right]g'=\g_1\left[\begin{array} {c} \om \\ 0
\end{array} \right]\]
where $\g_1=\g\left[\begin{array}{ll}  g&0\\
0& I_{n-m}
\end{array}\right]$ and $\om \in \frM_{k,m}$. If $\det(\g_1)=1$, we are
done. If $\det(\g_1)=-1$, one should
 replace  $\g_1 $ by  $\g_1 e, \; e=\left[\begin{array}{ll}  -1&0\\
0& I_{n-1}
\end{array}\right]$ and change $\om$ appropriately. The representation (\ref{form2})
 follows from (\ref{form1}).
\end{proof}

\begin{corollary} The integral (\ref{z-mes}) can be written as
\be\label{zzm} (\z_k,f)=\intl_{\text{\rm rank} (x) \le \, k} f(x)
d\nu_k (x)=\intl_{\text{\rm rank} (x) = k} f(x) d\nu_k (x).\ee
\end{corollary}
 Indeed, the first equality follows from Lemma \ref{nu} (i). The
 second one is clear from the observation that if $\rank (x) \le
 k-1$ then by (\ref{form2}), $x=v\om, \; v \in V_{n,k-1},
\; \om \in \frM_{k-1,m}$. The set of all such pairs $(v,\om)$ has
measure $0$ in $V_{n,k} \times \frM_{k,m}$.

\end{remark}

\section{Convolutions  with zeta distributions and Riesz potentials}

\subsection{Definitions} For $\a$ belonging to the Wallah set
(\ref{wal}), let us consider
 the convolution operator  $\z_\a \ast $ defined by
\be\label{al} (\z_\a \ast f)(x) =
\frac{1}{\gm(\a/2)}\intl_{\frM_{n,m}} f(x-y)|y|_m^{\a -n} dy\ee if
$Re \, \a>m-1$, and \be\label{alk}(\z_k \ast f)(x)=
c_1\intl_{V_{n,k}} dv\intl_{\frM_{k,m}} f(x-v \om)d\om, \ee
\[ c_1= 2^{-k} \, \pi^{(nm-km-nk)/2} \,
\Gam_k(n/2)/\Gam_m (n/2). \] if $\a=k, \; k=1,2, \ldots , n$; cf.
Corollary \ref{nf}. Note that for $m-1<k\le n$, both
representations are applicable to $\z_k \ast f$.

Another important normalization of the zeta distribution and the
corresponding convolutions (\ref{al}) and (\ref{alk}) leads to the
Riesz distributions and Riesz potentials. We observe that the
functional equation
 (\ref{z-f}) for the zeta function can be written in the form
 \be\label{z-f1}  \frac{1}{\g_{n,m}(\a)}\Z(f, \a-n)=(2\pi)^{-nm}
\Z(\F f, -\a)\; , \ee \be\label{gam} \gam_{n,m} (\a)=\frac{2^{\a
m} \, \pi^{nm/2}\, \Gam_m (\a/2)}{\Gam_m ((n-\a)/2)},\qquad \a\neq
n-m+1, \,  n-m+2, \ldots . \ee Now we have  excluded the values
$\a= n-m+1, \, n-m+2, \ldots $, because these are poles of the
gamma function $\Gam_m ((n-\a)/2)$ sitting in the numerator in
(\ref{gam}) (we recall that $m\geq 2$). The corresponding {\it
Riesz distribution} $h_\a$ is defined by
 \be\label{rd} (h_\a, f)= \frac{1}{\g_{n,m}(\a)}\Z(f, \a-n)= a.c. \; \frac{1}{\gam_{n,m} (\a)}
  \intl_{\Ma} |x|^{\a-n}_m f(x) dx,
  \ee
 where $f\in\S(\Ma)$. For $Re \, \a>m-1$, the distribution $h_\a$ is
 regular and agrees with the ordinary function $h_\a(x)=|x|^{\a-n}_m
 /\gam_{n,m}(\a)$.

 The {\it Riesz potential} of a function $f\in\S(\Ma)$
 is defined by \be\label{rie-z} (I^\a
f)(x)=(h_\a, f_x)=\frac{1}{\g_{n,m}(\a)}\Z(f_x, \a-n), \qquad
f_x(\cdot)= f(x-\cdot).\ee For $Re \, \a>m-1$, $ \a\neq n-m+1, \,
n-m+2, \ldots \;$, (\ref{rie-z}) is represented in the classical
form by the absolutely convergent integral
 \be\label{rie} (I^\a f)(x)=\frac{1}{\gam_{n,m} (\a)} \intl_{\Ma}
f(x-y) |y|^{\a-n}_m dy.\ee This integral operator is well known in
the rank-one case $m=1$.

 Riesz potentials of integral order deserve special mentioning.
The following representations are inherited from those  for the
normalized zeta function; see Theorem \ref{tzk}.
\begin{theorem}\label{ldes}  Let $f \in \S(\Ma)$. Suppose that
$\a=k$ is a positive integer. If $k \neq n-m+1, \, n-m+2, \ldots
$, then

 \bea(I^k f)(x)\label{des}&=&\g_1\intl_{SO(n)} d\gam \intl_{\frM_{k,m}}
 f \left (x-\gam \left[\begin{array} {c} \om \\ 0
\end{array} \right]  \right ) \, d\om, \\
\label{des1}&=& \g_2\intl_{V_{n,k}} dv\intl_{\frM_{k,m}} f(x-v
\om)d\om, \eea  where \bea\label{c1} \qquad \quad \g_1&=&2^{-km}
\,\pi^{-km/2}
\, \Gam_m \left(\frac{n-k}{2}\right) / \Gam_m \left(\frac{n}{2}\right),  \\
\label{c2}
\g_2&=& 2^{-k(m+1)} \, \pi^{-k(m+n)/2} \, \Gam_k
\left(\frac{n-m}{2}\right). \eea
\end{theorem}

The constant $\g_2$ above  was evaluated by making use of (\ref{2.5.2}).

The following theorem resumes basic properties of Riesz
distributions and Riesz potentials.
\begin{theorem}
Let $f\in\S(\Ma)$,  $ \a\in \bbc$,  $ \a\neq n-m+1, \, n-m+2,
\ldots \quad$.

\noindent {\rm (i)}  The Fourier transform of the Riesz
 distribution $h_\a$ is evaluated by the formula $(\F
 h_\a)(y)=|y|_m^{-\a}$, the precise meaning of which is
\be\label{fou} (h_\a, f)=(2\pi)^{-nm}(|y|_m^{-\a},(\F
f)(y))=(2\pi)^{-nm} \Z(\F f, -\a). \ee

\noindent {\rm (ii)} If $k=0,1,\dots$, and $\Del$ is the
Cayley-Laplace operator, then \be \label{Dkh}(-1)^{mk}\Del^k
h_{\a+2k}= h_{\a}, \quad i.e. \quad (-1)^{mk}( h_{\a+2k}, \Del^k
f)=( h_{\a}, f),\ee and, therefore, \be \label{Dkf}(I^{-2k}
f)(x)=(-1)^{mk}(\Del^k f)(x). \ee In particular, \be\label{I0}
(I^{0} f)(x)=f(x). \ee
\end{theorem}
\begin{proof}
 (i) follows immediately from the definition (\ref{rd}) and the
functional equation (\ref{z-f1} ). To prove (\ref{Dkh}), for
sufficiently large $\a$ we have \bea \Del^k
h_{\a+2k}(x)&=&\frac{1}{\g_{n,m}(\a+2k)}\Del^k |x|^{\a+2k-n}_m
\nonumber \\ &=&\frac{B_k(\a)}{\g_{n,m}(\a+2k)}\Del^k
|x|^{\a-n}_m\nonumber \\ &=& ch_{\a}(x),\nonumber \eea where by
(\ref{bka}) and (\ref{Poh}),
$$
c=\frac{B_k(\a)\g_{n,m}(\a)}{\g_{n,m}(\a+2k)}=\frac{B_k(\a)\Gam_m
(\a/2)\Gam_m ((n-\a)/2-k)}{4^{mk}\Gam_m (\a/2+k)\Gam_m
((n-\a)/2)}=1.
$$
For all admissible $\a\in\bbc$, (\ref{Dkh}) follows by analytic
continuation. The equality (\ref{Dkf}) is a consequence of
(\ref{fou}) and (\ref{Dkh}). Indeed, by (\ref{rie-z}),

\bea\nonumber(I^{-2k} f)(x)&=&(h_{-2k}, f_x)\stackrel{\rm
(\ref{Dkh})}{=}(-1)^{mk}(h_0, \;\Del^k
f_x)\\[14pt]\nonumber &\stackrel{\rm
(\ref{fou})}{=}&(-1)^{mk}(2\pi)^{-nm}\Z(\F(\Del^k f_x),\;
0)\\[14pt]\nonumber&=&(-1)^{mk}(\Del^k f_x)(0)=(-1)^{mk}(\Del^k
f)(x).\eea
\end{proof}

\subsection{Riesz potentials and  heat kernels} It is convenient to
study Riesz potentials by making use of the heat kernels and the
corresponding Gauss-Weierstrass integrals. In the rank-one case
$m=1$ this approach is described in \cite{Ru1}, Section 16. In the
higher rank case it was implicitly indicated in \cite{Cl} and
\cite{FK}. The key idea is to represent the Riesz potential as a
 lower-dimensional fractional integral of the
corresponding Gauss-Weierstrass integral which is  easy to handle.

\begin{definition} For $x \in \Ma, \;  n \ge m$, and $t \in \p$, we
define the (generalized) heat kernel $h_t (x)$ by the formula
\be\label{heat} h_t (x)=(4\pi)^{-nm/2}|t|^{-n/2} \exp (-\tr
(t^{-1} x'x)/4) \ee where $|t|=\det (t)$. The corresponding
Gauss-Weierstrass integral $(W_t f)(x)$ of a function $f(x)$ on
$\Ma$ is defined by  \be\label{ga} (W_t f)(x)=\intl_{\Ma} h_t
(x-y) f(y) dy=\intl_{\Ma} h_{I_m} (y) f(x-yt^{1/2}) \, dy. \ee
\end{definition}
\begin{lemma}{}\hfil

\noindent {\rm (i)} For each $ \; t \in \p$, \be\label{ed}
\intl_{\Ma} h_t (x) dx =1.\ee

\noindent {\rm (ii)} The Fourier transform of $ \; h_t (x)$ has
the form \be\label{ft}(\F h_t)(y)= \exp (-\tr (ty'y), \ee which
implies the semi-group property  \be\label{cnv} h_t \ast
h_\tau=h_{t+\tau}, \qquad t, \tau \in \p. \ee

\noindent {\rm (iii)} If $f \in L^p(\Ma), \; 1\le p \le \infty$,
then \be\label{gw} ||W_t f||_p \le ||f||_p \, , \qquad \quad W_t
W_\tau f=W_{t+\tau}f, \ee and \be\label{lim}\lim\limits_{t \to
0}(W_t f)(x)=f(x) \ee
 in the $L^p$-norm. If $f$ is  a continuous function vanishing at infinity, then
 (\ref{lim}) holds  in the $\sup$-norm.
\end{lemma}
\begin{proof} To prove (i), we pass to polar coordinates. Owing to
 (\ref{2.4}) and  (\ref{2.16}), we obtain \bea \intl_{\Ma} h_t (x) dx
&=&\frac{2^{-m} \, \sig_{n,m}}{(4\pi)^{nm/2} \,
|t|^{n/2}}\intl_{\p} \exp (-\tr (t^{-1}r)/4) \, |r|^{n/2 -d} dr
\nonumber \\ &=&\frac{2^{-m} \, \sig_{n,m}  \, \Gam_m (n/2)
}{(4\pi)^{nm/2} \, |t|^{n/2}} \, |t^{-1}/4|^{-n/2}= 1.\nonumber
\eea The formula (\ref{ft}) is well known (see, e.g., [Herz]), and
the proof is elementary. Indeed, by changing variable $x \to
2xt^{1/2}$ and setting $\eta=2yt^{1/2}$, we have
\[ (\F h_t)(y)=\pi^{-nm/2}\intl_{\Ma} \exp (-\tr (tx'x)) \, \exp
(i\tr (\eta'x)) \,  dx.\] This splits into a product of the
 one-dimensional Fourier transforms of Gaussian functions. The
 statement (iii) follows from (\ref{ed}) and (\ref{cnv}). Indeed, the relations in
 (\ref{gw}) are clear. Furthermore,
 $$ (W_t f)(x)- f(x)=\intl_{\Ma} h_{I_m} (y) [f(x-yt^{1/2}) -
 f(x)] \, dy.$$
Hence, by the generalized Minkowski inequality,
$$ ||W_t f- f||_p \le \intl_{\Ma} h_{I_m} (y) \, ||f(\cdot -yt^{1/2}) -
 f(\cdot)||_p  \, dy.$$
If $t$ tends to the $0$-matrix, then $yt^{1/2} \to 0$ in
$\bbr^{nm}$. Since the integrand above does not exceed $2||f||_p
\, h_{I_m} (y)$, we can pass to the limit under the sign of
integration, and the desired result follows. For continuous
functions vanishing at infinity, the argument is similar.
\end{proof}

\begin{remark}\label{hl} A challenging open problem is whether
 (\ref{lim}) holds  for almost all $x \in \Ma$. This fact is well known in the case $m=1$
[St2], [SW]. It  follows from the  estimate \be\label{max}
\sup_{t>0}|(W_t f)(x)| \le (M^*f)(x)\ee where $(M^*f)(x)$ is the
Hardy-Littlewood maximal function. It would be desirable to
extend this theory to the matrix case when the positive parameter
$t$ is replaced by a positive definite matrix.
\end{remark}

Let us consider the G{\aa}rding-Gindikin fractional integrals on
$\p$ defined by  \be\label{gg} (I_{-}^\lam g)(t) = \frac
{1}{\Gam_m (\lam)} \intl_t^\infty g(\t)|\tau -t |^{\lam-d} d\t,
\ee where $d=(m+1)/2, \; Re \, \lam > d-1$, and integration is
performed over all $\t \in \p$ so that $\t -t \in \p$ . The
following theorem establishes connection between the Riesz
potentials, the Gauss-Weierstrass integrals, and the
G{\aa}rding-Gindikin fractional integrals.
\begin{theorem}\label{rrg} Let $ \; m-1<Re \, \a<n-m+1$, $d=(m+1)/2$. Then
\be\label{rg} (I^\a f)(x) = \frac {1}{\Gam_m(\a/2)} \intl_{\p}
|t|^{\a/2-d}(W_t f)(x) dt, \ee \be\label{rgg} W_t [I^\a f](x)
=I_{-}^{\a/2}[(W_{(\cdot)} f)(x)](t)\ee provided that integrals on
either side of the corresponding equality exist in the Lebesgue
sense.
\end{theorem}
\begin{proof} The right-hand side of (\ref{rg}) transforms as
follows:
\[\text{\rm r.h.s.}=\frac{(4\pi)^{-nm/2}}{\Gam_m(\a/2)}  \intl_{\p}
|t|^{(\a-n)/2-d} dt\intl_{\Ma} f(y) \exp (-\tr (t^{-1}z)) \, dy\]
where $z=(x-y)'(x-y)/4$. Now we change the order of integration
and replace $t^{-1}$ by $t$. Since
\[ \intl_{\p} |t|^{(n-\a)/2-d} \exp (-\tr (tz)) \, dt=\Gam_m
\left(\frac{n-\a}{2}\right) |z|^{(\a-n)/2}, \qquad |z| \neq 0,\]
 (\ref{rg}) follows. The validity of (\ref{rgg}) is a simple
 consequence of the semi-group property of the Gauss-Weierstrass
 integral. Indeed, \bea W_t [I^\a f](x)
&=&\frac {1}{\Gam_m(\a/2)} \intl_{\p} |\t|^{\a/2-d}(W_t W_\t
f)(x)\,
d\t \nonumber \\
&=&\frac {1}{\Gam_m(\a/2)} \intl_{\p} |\t|^{\a/2-d}(W_{t+\t}
f)(x)\,
d\t \nonumber \\
&=&I_{-}^{\a/2}[(W_{(\cdot)} f)(x)](t).\nonumber \eea
\end{proof}

Theorem \ref{rrg} has a number of remarkable consequences.
Firstly, it enables us to invert the Riesz potential by inverting
 $I_{-}^{\a/2}$ and
applying the  approximation property $\lim\limits_{t \to 0}W_t
f=f$. We shall study this question  in subsequent publications.
Secondly, Theorem \ref{rrg} provides a simple proof of   the
semigroup property of the Riesz potentials under mild assumptions
for $f$.
\begin{corollary} Let \be\label{cond} Re \, \a>m-1, \quad Re \, \b>m-1, \quad Re \,
(\a+\b)<n-m+1. \ee Then \be\label{sg} (I^\a I^\b f)(x)= (I^{\a+\b}
f)(x) \ee provided that the integral $(I^{\a+\b} f)(x)$ absolutely
converges.
\end{corollary}
\begin{proof} Applying Theorem
\ref{rrg} and changing the order of integration, we obtain \bea
I^\a I^\b f&=&\frac {1}{\Gam_m(\a/2)} \intl_{\p} |t|^{\a/2-d} \,
W_t I^\b f
\,  dt \nonumber \\
&=&\frac {1}{\Gam_m(\a/2)} \intl_{\p} |t|^{\a/2-d}
\,(I_{-}^{\b/2}W_{(\cdot)} f)(t)\,  dt \nonumber \\
&=&\frac {1}{\Gam_m(\a/2) \, \Gam_m(\b/2)} \intl_{\p} W_\t f \,
d\t \intl_0^\t |t|^{\a/2-d}|\t -t|^{\b/2-d}
\,  dt \nonumber \\
&=&\frac {1}{\Gam_m((\a+\b)/2)} \intl_{\p} |\t|^{(\a+\b)/2-d}
\,W_{\t} f \, d\t=I^{\a+\b} f.\nonumber \eea
\end{proof}
\begin{remark} For $m=1$, the conditions (\ref{cond}) have a
well-known form
\[ Re \, \a>0, \quad Re \, \b>0, \quad Re \,
(\a+\b)<n. \] Note that in the higher rank case, (\ref{cond}) is
possible only if $3m<n+3$. It is also worth noting that the
``usual" way to prove the semi-group property (\ref{sg}) for
``rough" functions $f$ via the Fourier transform formula and
evaluation of the corresponding beta integral  (cf. [SKM], [St2])
becomes much more difficult in the higher rank case.
\end{remark}
\begin{remark}\label{arc}
Since the ``symbol" (i.e., the Fourier transform of the kernel) of
the Riesz potential $I^\a f$ is $|y|_m^{-\a}$, then, for
sufficiently good functions $f$, (\ref{sg}) holds for arbitrary
complex $\a $ and $\b$. We will be concerned with this topic in
Section 6.
\end{remark}

\subsection{$L^p$-convolutions} Another useful application of Theorem
\ref{rrg} is the following.
 By Lemma \ref{lacz} and Theorem \ref{tzk},  convolutions  with zeta
 distributions and  Riesz potentials  are well defined on test functions belonging to
the Schwartz space. Are they meaningful for $f \in L^p \, $?  Let
us study this question.

\begin{theorem}\label{ap} For $f \in L^p(\Ma)$, the Riesz potential $(I^\a
f)(x)$ absolutely converges almost everywhere on $\Ma$ provided
\be\label{sh} Re \, \a > m-1, \qquad 1 \le p <\frac{n}{Re \, \a
+m-1}. \ee
\end{theorem} \begin{proof} Note that (\ref{sh}) is possible if
$m \le [(n+1)/2]$ where $[(n+1)/2]$ is the integral part of
$(n+1)/2$. The conditions (\ref{sh}) are well known if $m=1$ [St1]
 when they have the form $
 1 \le p <n/Re \, \a$. To prove our theorem,  it suffices
to consider $f \ge 0$ and  justify the following inequality
\be\label{ine} I \equiv \intl_{\Ma} \exp (-\tr (x'x)) \, (I^\a
f)(x) \, dx \le c \, ||f||_p \, , \quad c=\const.\ee  Regarding
$I$ as the Gauss-Weierstrass
 integral $(W_{I_m}[I^\a f])(0)$ owing  to (\ref{rgg}), we have
 \bea I&{}\simeq {}&\intl_{I_m}^\infty |\t -I_m|^{\a/2 -d} (W_\t f)(0) \, d\t
  \nonumber \\
 &{}\simeq {}&\intl_{I_m}^\infty |\t -I_m|^{\a/2 -d} \, |\t|^{-n/2}  \,d\t
 \intl_{\Ma} \exp (-\tr (\t^{-1}x'x)/4) \, f(x) \, dx,\nonumber \eea
 $d=(m+1)/2$. By  H\"older's inequality, $I{}\lesssim{}A||f||_p \; $
 where \[
 A^{p'}=\intl_{I_m}^\infty |\t -I_m|^{\a/2 -d} \, |\t|^{-n/2}
\Big ( \,\intl_{\Ma} \! \exp (-p'\tr (\t^{-1}x'x)/4)\, dx \Big
)^{1/p'}d\t, \] $1/p +1/p' =1$. This integral can be easily
estimated. Indeed, by setting  $x=y\t^{1/2}, \; dx=|\t|^{n/2}dy$,
we obtain \bea A^{p'}
&{}\simeq {}&\intl_{I_m}^\infty |\t -I_m|^{\a/2 -d} \, |\t|^{-n/2p} \,d\t \qquad (\t^{-1}=r) \nonumber \\
&{}\simeq {}&\intl_0^{I_m} |I_m -r|^{\a/2 -d} \, |r|^{(n/p
-\a)/2-d} \,dr. \nonumber \eea The last  integral is finite if
$p$ obeys (\ref{sh}).
\end{proof}

More precise information can be obtained for convolutions with
zeta distributions (and therefore, for Riesz potentials) of
integral order. For $ k=1,2, \dots \, n-1$,
 we denote \be\label{den} v_{n-k} \! = \! \left[\begin{array} {c} 0 \\
I_{n-k}
\end{array} \right] \! \in  \! V_{n,n-k}, \quad \Pr_{\bbr^{n-k}} \! = \! v_{n-k}v'_{n-k} \! = \!  \!
\left[\begin{array}{ll}  0&0\\
0&I_{n-k}
\end{array}\right],\ee
\be\label{rad} \tilde f (x)=\intl_{SO(n)} f(\gam x) d\gam.\ee
\begin{lemma}\label{zz} Let $n \ge m \ge 1; \;  k=1,2, \ldots , n$. For any
$\lam>k+m-1$, \be\label{rav} \intl_{\Ma} \frac{|(\z_k \ast
f)(x)|}{|I_m + x'x|^{\lam/2}} \, dx= c_\lam\intl_{\Ma} \frac{
\tilde f (x)}{|I_m + x' \Pr_{\bbr^{n-k}} x|^{(\lam -k)/2}} \,
dx,\ee \be\label{con}c_\lam= \frac{\pi^{nm/2} \, \Gam_m ((\lam
-k)/2)} {\Gam_m (n/2) \, \Gam_m (\lam /2)}. \ee
\end{lemma}
\begin{proof} Suppose first that $k<n$, and denote the left-hand side of
(\ref{rav}) by $I(f)$. By (\ref{alk}), \bea
 I(f)&=&c_1\intl_{V_{n,k}}  \,dv\intl_{\frM_{k,m}}d\om \intl_{\Ma}
f  (x-v\om)\frac{dx}{|I_m + x'x|^{\lam/2}} \nonumber \\
&{}&(\text{\rm set} \quad v=\g v_0, \;  v_0=\left[\begin{array} {c} I_k \\
0
\end{array} \right], \;  x=\gam y, \quad \g \in SO(n)) \nonumber \\
&=&c_1 \, \sig_{n,k}\intl_{\frM_{k,m}}d\om \intl_{\Ma}
\tilde f  (x-v_0\om)\frac{dx}{|I_m + y'y|^{\lam/2}} \nonumber \\
&=&c_1 \, \sig_{n,k}\intl_{ \frM_{n-k,m}} db
\intl_{\frM_{k,m}}\tilde f \left ( \left[\begin{array} {c} \om \\
b
\end{array} \right]  \right )d\om\intl_{\frM_{k,m}}
\frac{da}{|I_m + b'b +a'a|^{\lam/2}}. \nonumber \eea For
$\lam>k+m-1$, the last integral can be explicitly evaluated by the
formula \be\label{for}\intl_{\frM_{k,m}}\frac{da}{|q
+a'a|^{\lam/2}}= \frac{\pi^{km/2} \, \Gam_m ((\lam -k)/2)} {
\Gam_m (\lam /2)} |q|^{(k-\lam)/2}, \qquad q \in \p \ee (see
(A.3)). This gives (\ref{rav}). For $k=n$, the proof is
simpler and based on (\ref{for}). In this case $ \Pr_{\bbr^{n-k}}$
should  be replaced by the zero matrix.
\end{proof}

\begin{corollary} If $n \ge m \ge 1$,  $  k=1,2, \ldots , n$, and
$\lam>k+m-1$, then \be\label{rav1} \intl_{\Ma} \frac{|(\z_k \ast
f)(x)|}{|I_m + x'x|^{\lam/2}} \, dx \le c_\lam ||f||_1.\ee If
$k=n$ and $f \ge 0$,  then (\ref{rav1}) is a strict equality.
\end{corollary}

Let us extend (\ref{rav1}) to $f \in L^p$ with $p \ge 1$. Now we
have to impose extra restrictions because of the projection
operator $\Pr_{\bbr^{n-k}}$ on the right-hand side of
 (\ref{rav}).
\begin{theorem}\label{lp} Let $n>m; \;   k=1,2, \ldots , n-m$. If $f \in L^p$ and
$$
\lam >k+\max \left (m-1, \frac{n+m-1}{p'}\right ), \qquad
\frac{1}{p}+ \frac{1}{p'}=1, $$ then \be\label{rav2} \intl_{\Ma}
\frac{|(\z_k \ast f)(x)|}{|I_m + x'x|^{\lam/2}} \, dx \le c
||f||_p, \qquad c=c(\lam, p, n,k,m),\ee provided \be\label{pp}1
\le p<\frac{n}{k+m-1}.\ee
\end{theorem}
\begin{proof} Note that (\ref{pp})  agrees with (\ref{sh})
 in Theorem \ref{ap}. It suffices to prove the statement for
non-negative radial functions $f(x)\equiv f_0(x'x)$.  We remind
the notation $d=(m+1)/2$.  Denote by $I(f)$ the left-hand side of
(\ref{rav2}) and make use of
 (\ref{rav}). Splitting the integral in the right-hand side of (\ref{rav})
 in $ \int_{\frM_{n-k,m}} \times  \int_{\frM_{k,m}}$, and passing in
$\int_{\frM_{n-k,m}}$ to polar coordinates, we have \bea
I(f)&{}\simeq {}&\intl_{\p}\frac{|r|^{(n-k)/2 -d}} {|I_m
+r|^{(\lam -k)/2}} \, dr \intl_{\frM_{k,m}} f_0(\om'\om +r) \,
d\om
\nonumber \\
&=&\intl_{\frM_{k,m}} d\om \intl_{\om'\om}^\infty f_0(s) \,
 \frac{|s-\om'\om |^{(n-k)/2 -d} }{|I_m +s-\om'\om|^{(\lam -k)/2}} \, ds \nonumber \\
&=&\intl_{\p} f_0(s) \, ds \intl_{\om'\om <s} \frac{|s-\om'\om |^{(n-k)/2 -d} }{|I_m +s-\om'\om|^{(\lam -k)/2}} \, d\om\nonumber \\
&&(\text{\rm set} \quad \om=\eta s^{1/2}, \quad \eta \in
\frM_{k,m}, \quad
d\om=|s|^{k/2} \, d\eta ) \nonumber \\
&=&\intl_{\p} f_0(s) |s|^{n/2 -d}ds \intl_{\eta'\eta <I_m}
\frac{|I_m-\eta'\eta |^{(n-k)/2 -d} } {|I_m
+s^{1/2}(I_m-\eta'\eta)s^{1/2}|^{(\lam -k)/2}}\, d\eta. \nonumber
\eea By taking into account that \be\label{nor} ||f||_p{}\simeq {}
\left ( \, \intl_{\p} | f_0(s)|^p  |s|^{n/2 -d}ds\right )^{1/p},
\ee owing to the H\"older and the generalized Minkowski
inequality, we obtain $I(f){}\lesssim{}c ||f||_p$ where
\[ c=\intl_{\eta'\eta <I_m}
|I_m-\eta'\eta |^{(n-k)/2 -d} A_\mu^{1/p'} (I_m-\eta'\eta) \,
d\eta, \]
\[ A_\mu (r)=\intl_{\p}\frac {|s|^{n/2 -d}}{|I_m +s^{1/2}rs^{1/2}|^{\mu/2}} \, ds,
\quad \mu=(\lam -k)p', \quad  r=I_m-\eta'\eta.\] Since $|I_m
+s^{1/2}rs^{1/2}|=|I_m +r^{1/2}sr^{1/2}|$, by changing variable
$s=r^{-1/2}\rho r^{-1/2}$, we have
\[ A_\mu (r)=|r|^{-n/2}\intl_{\p}\frac {|\rho|^{n/2 -d}}
{|I_m +\rho |^{\mu/2}} \, d\rho.\] The last integral is finite for
$\mu>n+m-1$ or $(\lam -k)p'>n+m-1$ (moreover, it can be explicitly
evaluated
 by (A.2)). Thus for $\lam
>k+(n+m-1)/p'$ we have \bea c&{}\simeq {}&\intl_{\eta'\eta
<I_m}|I_m-\eta'\eta |^{(n-k)/2- n/2p' -d} d\eta \nonumber \\
&{}\simeq {}&\intl_0^{I_m}|I_m-r|^{(n/p-k)/2-d}|r|^{-1/2} \, dr<
\infty\nonumber \eea provided $n/p -k>m-1$, i.e., $1 \le p
<n/(k+m-1)$. This completes the proof.
\end{proof}
\begin{remark}\label{rpo} We do not know whether the conditions (\ref{pp}) and
(\ref{sh}) are sharp. Moreover, it would be interesting to study
boundedness properties of the Riesz potential operator $I^\a$ on
functions $f \in L^p(\Ma)$ in different norms.
\end{remark}

\section {Radon transforms}

\setcounter{equation}{0}

\subsection {Definitions and basic properties}
  Let $k,n$, and $m$ be
positive integers, $0<k<n$, $\vnk$ be the Stiefel manifold of
orthonormal $(n-k)$-frames in $\bbr^n$. For $\; \xi\in\vnk$, and
$t\in\Mt$, the linear manifold \be\label{pl} \tau=
\tau(\xi,t)=\{x\in\Ma:\eq\} \ee
 will be called {\it a matrix $k$-plane} in $\Ma$. We denote by  $\frT$  the
  variety of all such
 planes. For $m=1$, $\frT$ is an affine Grassmann manifold of $k$-dimensional planes $\bbr^n$.
The matrix $k$-plane Radon transform $\hat f(\tau)$ of a function
$f(x)$ on $\Ma$ assigns to $f$ a collection of integrals of $f$
over all matrix planes $ \tau \in \frT$.
 Namely, \[ \hat f (\tau)=\int_{x \in \tau} f(x), \qquad \tau \in \frT.\]
 The  dual Radon transform
$\check \vp (x)$ is a mean value of a function $\vp(\tau)$ on
$\frT$ over all matrix planes $\tau$ through $x$:
$$
\df=\intl_{\tau\ni x}\vp(\tau), \qquad x\in\Ma.
$$
Precise meaning of these integrals is the following. Let $\tau$ be
the plane (\ref{pl}), $ \xi_0=\left[\begin{array} {c}  0 \\
I_{n-k} \end{array} \right] \in \vnk$, and $ g_\xi \in SO(n)$ be a
rotation satisfying $ g_\xi\xi_0=\xi$. Then \be\label{4.9} \hat f
(\tau) \equiv \rf=\intl_{\Mkm} f\left(g_\xi \left[\begin{array}
{c} \om \\t
\end{array} \right]\right)d\om,\ee and
\be\label{4.2}  \df=\frac{1}{\sigk}\intl_{\vnk}
\varphi(\xi,\xi'x)d\xi=\intl_{SO(n)}\varphi(\g\xi_0,\xi_0'\g
'x)d\g. \ee
\begin{lemma}\label{dur} {}\hfill

\noindent {\rm (i)} The duality relation \be\label{4.3}
\intl_{\Ma} f(x)\df dx=\frac{1}{\sigk}
\intl_{\vnk}d\xi\intl_{\Mt}\fc\rf dt \ee holds provided that
either side of (\ref{4.3}) is finite for $f$ and $\vp$ replaced by
$|f|$ and $|\vp|$, respectively.

\noindent {\rm (ii)} If $f\in L^1(\Ma)$, then the Radon transform
$\rf$ exists for all $\xi\in\vnk$ and almost all $t\in\Mt$.
Furthermore, \be \intl_{\Mt}\rf \, dt=\intl_{\Ma} f(x) \, dx,
\qquad \forall  \, \xi\in\vnk.\ee
\end{lemma}
\begin{proof}
By (\ref{4.9}), the right-hand side of (\ref{4.3}) has the form
\be\label{rs} \frac{1}{\sigk} \intl_{\vnk}d\xi\intl_{\Mt}\fc \, dt
\intl_{\Mkm} f\left(g_\xi \left[\begin{array} {c} \om \\t
\end{array} \right]\right)d\om.\ee
Changing variables $x=g_\xi \left[\begin{array} {c} \om \\t
\end{array} \right]$, we have
\[ \xi'x=(g_\xi \xi_0)'g_\xi \left[\begin{array} {c} \om \\t
\end{array} \right]=\xi_0'\left[\begin{array} {c} \om \\t
\end{array} \right]=t.\]
Hence,  by the Fubini theorem,  (\ref{rs}) reads \[
\frac{1}{\sigk} \intl_{\vnk}d\xi\intl_{\Ma} \varphi(\xi,\xi'x)
f(x) \, dx=\intl_{\Ma} \df  f(x) \, dx.\] The statement (ii) is a
consequence of the Fubini theorem: \bea \intl_{\Mt}\rf
dt&=&\intl_{\Mt} dt \intl_{\Mkm} f\left(g_\xi \left[\begin{array}
{c} \om \\t
\end{array} \right]\right)d\om \nonumber \\
 &=&\intl_{\Ma} f(g_\xi x) \,
dx=\intl_{\Ma} f(x) \, dx. \nonumber \eea
\end{proof}
 The following
properties  of the Radon transform can be easily checked.
\begin{lemma}\label{l4.3} Suppose that the Radon transform \[ f
(x) \longrightarrow \rf, \qquad x \in \Ma, \quad (\xi, t) \in \cd,
\]
 exists (at least almost everywhere). Then
 \be\label{ev} \hat f(\xi\theta ',\theta t)=\rf,  \qquad \forall \theta\in O(n-k). \ee
 If $ \; g(x)=\gamma x\beta +y \; $ where $ \;
\gamma\in O(n), \quad  \b\in O(m), \quad  y\in\Ma$, then
\be\label{4.23} (f \circ g)^\wedge (\xi, t)= \hat
f(\gamma\xi,t\beta +\xi'\gamma 'y). \ee In particular, if
 $ \; f_y(x)=f(x+y)$, then \be
\label{4.4}
 \hat f_y(\xi ,t)= \hat f(\xi , \xi 'y+t).  \ee
\end{lemma}

The equality (\ref{ev}) is a matrix analog of the ``evenness
property" of the classical Radon transform, cf.  [Hel].

It is known [OR1] that the Radon transform  $f \to \hat f$ is
injective on the Schwartz space $\S=\S(\Ma)$ if and only if $k \le
n-m$. The classical problem is how to reconstruct a function $f$
from the integrals $\hat f (\tau)$,  $\t \in \frT$. One of the
ways  to do this lies via Riesz potentials $I^k f$ of integral
order. By Theorem \ref{ldes}, for  $1\le k \le  n-m$ and $f \in
\S$ we have \be\label{any} (I^k f)(x) = \g_1 \intl_{SO(n)} d\gam
 \intl_{\frM_{k,m}} f \left (x-\gam
\left[\begin{array} {c} \om \\ 0
\end{array} \right]  \right ) \,
 d\om,\ee
\be\label{cr} \g_1=2^{-km} \,\pi^{-km/2} \, \Gam_m\Big (
\frac{n-k}{2}\Big ) / \Gam_m\Big ( \frac{n}{2}\Big ). \ee By
Theorem \ref{lp}, this expression is well defined a.e. on $\Ma$ as
an absolutely convergent integral for any $f \in L^p$ provided $1
\le p<n/(k+m-1)$. The following important statement establishes
connection between Riesz potentials and Radon transforms. It
generalizes the classical result of B. Fuglede; see \cite{Fu},
\cite {Hel} for $m=1$.
\begin{theorem}\label{main} Let $ 1 \le k \le n-m$.
 Then
 \be\label{fu} \g_1(\hat
f)^{\vee} (x) \! = \! (I^k f)(x),\ee  provided  the Riesz
potential  on the right-hand side absolutely converges, e.g., for
$f \in L^p$,  $1 \le p<n/(k+m-1)$.
\end{theorem}
\begin{proof}  Let $f_x
(y)=f(x+y)$. Then  (\ref{4.4}) yields $\hat f_x(\xi ,t)= \hat
f(\xi , \xi 'x+t)$. Owing to (\ref{4.2}) and (\ref{4.9}), by
changing the order of integration, we have \bea(\hat f)^{\vee}
(x)&=&\frac{1}{\sigk}\intl_{\vnk} \hat
f(\xi,\xi'x)d\xi=\frac{1}{\sigk}\intl_{\vnk} \hat f_x(\xi,0)d\xi
\nonumber \\
&=&\intl_{SO(n)} d\gam
 \intl_{\frM_{k,m}} f \left (x+\gam
\left[\begin{array} {c} \om \\ 0
\end{array} \right]  \right ) \,
 d\om.
\nonumber \eea  This coincides with $\g_1^{-1}(I^k f)(x)$.
\end{proof}
\begin{corollary} \label{corr}If $f \in L^p$,  $1 \le p<n/(k+m-1)$, then $\hat
f(\t)$ is finite for almost all $\t \in \frT$.
\end{corollary}
\begin{remark} The condition $1 \le p<n/(k+m-1)$ in Corollary \ref{corr} can be
replaced by the weaker one, namely, \be\label{llp} 1\leq
p<\frac{n+m-1}{k+m-1}; \ee see [OR2]. The proof of this fact
relies on completely different ideas related to representation of
the Radon transform of  radial functions by the
G{\aa}rding-Gindikin fractional integral (\ref{gg}). For  $m=1$, both conditions coincide and cannot be improved. 
We are not
sure that for $m >1$, (\ref{llp}) is sufficient for the existence
of the Riesz potential $I^k f$ (in contrast to the case $m=1$);
cf. Remark \ref{rpo}.
\end{remark}

\subsection{The inversion problem } Theorem \ref{main}  reduces the
inversion problem for the Radon transform to that for the Riesz
potentials (as in the rank-one case [Hel], [Ru2]). Below we show
how the unknown function $f$ can be recovered in the framework of
the theory of distributions.

Let us consider the Riesz distribution
 \be\label{dis} (h_\a, f)= a.c.\; \frac{1}{\gam_{n,m} (\a)}
 \intl_{\Ma} |x|^{\a-n}_m f(x) dx, \ee
where $f \in \S=\S(\Ma)$.  From the Fourier transform formula
$$(h_\a, f)=(2\pi)^{-nm}(|y|_m^{-\a},(\F f)(y)),$$
 it is evident that the Schwartz class $\S$
does not suit well enough
  because it is not invariant under multiplication by
$|y|_m^{-\a}$. To get around this difficulty, we choose another
space of test functions.
  Let $\Psi=\Psi(\Ma)$ be the collection
of all functions $\psi (y) \in \S(\Ma)$ vanishing on the manifold
\be\label{sets} V \!= \!\{y: \, y \in \Ma, \;  \rank (y)<m \} \! =
\!\{y: \,y \in \Ma, \; |y'y| \!= \!0 \}\ee with all derivatives
(the coincidence of both sets in (\ref{sets}) is clear because
$\rank (y)=\rank (y'y)$). The manifold $V$ is a cone in
$\bbr^{nm}$ with vertex $0$. The set $\Psi$ is a closed linear
subspace of $\S$. Therefore, it can be regarded as a linear
topological space with the induced topology of $\S$. Let $\Phi
=\Phi(\Ma)$ be the
 Fourier image of $\Psi$. Since the Fourier transform
$\F$ is an automorphism of $\S$ (i.e., a topological isomorphism
of $\S$ onto itself), then $\Phi$ is a closed linear subspace of
$\S$. Having been equipped with the induced topology of $\S$, the
space $\Phi$ becomes a linear topological space isomorphic to
$\Psi$ under the Fourier transform.

The spaces  $\Phi$ and $\Psi$ were introduced by V.I. Semyanistyi
[Se] in the rank-one case $m=1$.  They have proved to be very
useful in integral geometry, fractional calculus, and  the theory
of function spaces. Further generalizations and applications can
be found in [Li], [Ru1], [Sa], [SKM].

In our case the following characterization of the space $\Phi$ is
a consequence of a  more general result due to S.G. Samko [Sa].
\begin{theorem} The Schwartz function $\phi (x)$ on $\Ma$ belongs
to the space $\Phi$ if and only if it is orthogonal to all
polynomials  $p(x)$ on any hyperplane $\t$ in $\bbr^{nm}$ having
the form $\t=\{x: \tr (a'x)=c \}, \; a \in V$: \be \intl_{\t} p(x)
\phi (x)d\mu(x)=0,\ee $d\mu(x)$ being the induced Lebesgue measure
on $\t$.
\end{theorem}

Note that for $m=1$, the space $\Phi$ consists of Schwartz
functions which are orthogonal to all polynomials on $\bbr^{n}$.

 We denote by  $\Phi'$ the space of all linear
continuous functionals (generalized functions) on $\Phi$. It is
clear that two $\S'$-distributions that coincide in the
$\Phi'$-sense, differ from each other by an arbitrary
$\S'$-distribution with the Fourier transform supported by $V$.
Since for any complex $\a$, multiplication by $|y|_m^{-\a}$ is an
automorphism of $\Psi$, then, according to the general theory
[GSh2], $I^\a$ is an automorphism of $\Phi$, and we have \be
\label{fof}\F[I^\a \phi](y)=|y|_m^{-\a}\F[\phi](y), \qquad \phi
\in \Phi.\ee This
 gives \be\label{sgp} I^\a I^\b \phi =I^{\a+\b} \phi \quad
\forall \a, \b \in\bbc \ee (cf. Remark \ref{arc}), and \be
\label{ff}\F[I^\a f](y)=|y|_m^{-\a}\F[f](y)\ee for all
$\Phi'$-distributions $f$. The last  formula implies the following
\begin{theorem}\label{inv} The function
$ f \in L^p(\Ma), \; 1 \le p <n/(k+m-1)$, can be recovered from
the Radon transform $g=\hat f$  in the sense of
$\Phi'$-distributions by the formula \be (f,\phi)=\g_1(\check g,
I^{-k}\phi), \qquad \phi \in \Phi, \ee where
$$(I^{-k}\phi)(x)=(\F^{-1}|y|_m^{k}\F\phi)(x),$$
$\g_1$ being the constant (\ref{cr}).\end{theorem}
\begin{proof} We have \bea (f,\phi)&=&(2\pi)^{-nm} (\F[f],
\F[\phi])\nonumber \\ &=&(2\pi)^{-nm} (|y|_m^{-k}\F[f](y),
|y|_m^{k}\F[\phi](y)) \nonumber \\ &=&(2\pi)^{-nm} ((\F[I^k
f])(y), |y|_m^{k}\F[\phi](y))\nonumber \\ &=&(2\pi)^{-nm} \g_1
((\F[\check g])(y), |y|_m^{k}\F[\phi](y))\nonumber \\ &=& \g_1
(\check g, I^{-k} \phi).\nonumber \eea \end{proof}
\begin{remark}\label{clap} For $k$ even, the Riesz potential $I^kf$ can be
inverted (in the sense of $\Phi'$-distributions) by repeated
application  of the Cayley-Laplace operator $ \Del_m =\det
(\partial' \partial), \;  \partial=(\partial/\partial x_{i,j})$.
 This operator agrees with multiplication by
$(-1)^{m} |y|_m^2$   in the Fourier terms, and therefore,
$(-1)^{km} \Del_m^k I^k f= f$ in the $\Phi'$-sense.
\end{remark}

\begin{remark}\label{op3} It would be desirable to obtain {\it
pointwise
 inversion formulas} for $I^k
f $ and $\hat f$ (not in the $\Phi'$-sense).
 In the rank-one case such formulas can be found in [Ru1] and [Ru2].
\end{remark}

\section {Appendix}
The formulas in this section are not new. Since it is not so easy
to find simple proofs of them in the literature, we present such
proofs for convenience of the reader. We recall that $\p$ denotes
the cone of positive definite $m \times m$ matrices, $\cpm$ is the
closure of $\p$, $B_m$ is the beta function (\ref{2.6}),
$d=(m+1)/2$. The following formulas hold:
 \[\label{2.14} (A.1)\qquad  \qquad \qquad \intl_s^\infty
|r|^{-\gam}|r-s|^{\a-d}dr=|s|^{\a-\gam}B_m ( \a,\gam-\a),\qquad
\quad \]
$$ s\in\p, \quad Re\,\a>d-1, \quad Re\, (\gam -\a)> d-1;$$
\[\label{2.13}(A.2)  \qquad \qquad \intl_s^\infty
|I_m+r|^{-\gam}|r-s|^{\a-d}dr=|I_m+s|^{\a-\gam}B_m (
\a,\gam-\a),\qquad \quad \]
$$ \quad  s\in\cpm, \quad Re\,\a>d-1, \quad Re\, (\gam -\a)> d-1;$$
\[\label{2.15}(A.3) \qquad \qquad
\intl_{\Mkm}|b+y'y|^{-\lam/2}dy=\frac{\pi^{km/2}\gm((\lam-k)/2)}{\gm(\lam/2)}|b|^{(k-\lam)/2},
\qquad \quad  \]
$$ b\in\p, \quad Re\,\lam>k+m-1;$$
\[\label{2.15.1}(A.4)\qquad
   \intl_{\{y\in\Mkm:\;y'y<b\}}\!\!\!\!\!\!\!
|b-y'y|^{(\lam-k)/2-d
}dy=\frac{\pi^{km/2}\gm((\lam-k)/2)}{\gm(\lam/2)}|b|^{\lam/2-d},\qquad
\quad \]
 $$ b\in\p, \quad Re\,\lam>k+m-1.$$

 \noindent {\it Proof}.

 \noindent {\bf (A.1), (A.2).}
 By setting $r=q^{-1}$, $dr=|q|^{-m-1}dq$, one can write the left side of
(A.1) as \be |s|^{\a-d}\intl_0^{s^{-1}}
|q|^{\gam-\a-d}|s^{-1}-q|^{\a-d}dq\\\nonumber =|s|^{\a-\gam}B_m (
\a,\gam-\a), \ee and we are done. The equality  (A.2) follows from
(A.1) if we replace $s $ and $r$ by $I_m +s$ and $I_m +r$,
respectively.

 \noindent {\bf  (A.3), (A.4).}  By changing variable $y\to y
b^{1/2}$, we obtain
\[
\intl_{\Mkm}|b+y'y|^{-\lam/2}dy=|b|^{(k-\lam)/2}J_1,\] \[
\intl_{\{y\in\Mkm:\;y'y<b\}}
|b-y'y|^{(\lam-k)/2-d}dy=|b|^{\lam/2-d}J_2,\] where \bea\nonumber
J_1&=&\intl_{\Mkm}|I_m+y'y|^{-\lam/2}dy,\\
J_2&=&\intl_{\{y\in\Mkm:\;y'y<I_m\}}
|I_m-y'y|^{(\lam-k)/2-d}dy.\nonumber \eea Thus we have to show
that
\[ J_1=J_2=\frac{\pi^{km/2}\gm((\lam-k)/2)}{\gm(\lam/2)}.\]

 {\bf The case $k \ge m$.} We write both integrals in the polar coordinates according to
 Lemma \ref{l2.3}. For $J_1$ we have
\bea J_1&=&2^{-m}\sigma_{k,m}\intl_{\p} |r|^{k/2
-d}|I_m+r|^{-\lam/2}dr \nonumber \\
&=&2^{-m}\sigma_{k,m}B_m\left(\frac{k}{2},\frac{\lam-k}{2}\right)
\nonumber \eea (the second equality holds by (A.2) (with $s=0, \;
\a=k/2, \; \g=\lam/2 $). Similarly,
 \bea
J_2&=&2^{-m}\sigma_{k,m}\intl_{0}^{I_m}|r|^{k/2
-d}|I_m-r|^{(\lam-k)/2-d}dr \nonumber \\
&=&2^{-m}\sigma_{k,m}B_m\left(\frac{k}{2},\frac{\lam-k}{2}\right).\nonumber
\eea Now the result follows by (\ref{2.6}) and (\ref{2.16}).

{\bf The case $k < m$.} We replace $y$ by $y'$ and pass to the
polar coordinates. This yields
 \bea J_1&=&\intl_{\Mmk}|I_m+yy'|^{-\lam/2}dy \nonumber \\
&=&2^{-k}\intl_{\vmk}dv\intl_{\pk}
|I_m+vqv '|^{-\lam/2}|q|^{(m-k-1)/2}dq, \nonumber \\
&&(|I_m+vqv '|=|I_k+q|)\nonumber \\
&=& 2^{-k}\sigma_{m,k}
\intl_{\pk}|q|^{(m-k-1)/2}|I_k+q|^{-\lam/2}dq. \nonumber \eea By
(A.2) (with $s=0, \;  m=k, \; \g=\lam/2 $) and (\ref{2.5.2}), \bea
J_1&=&2^{-k}\sigma_{m,k}B_k\left(\frac{m}{2},\frac{\lam-m}{2}\right)\nonumber
\\&=& \frac{\pi^{km/2}\gk((\lam-m)/2)}{\gk(\lam/2)} \nonumber
\\&=&
\frac{\pi^{km/2}\gm((\lam-k)/2)}{\gm(\lam/2)}.\nonumber \eea
Similarly, \bea
J_2&=&\intl_{\{y\in\Mmk:\;yy'<I_m\}}|I_m-yy'|^{(\lam-k)/2-d}dy
\nonumber \\
&=&2^{-k}\intl_{\vmk}dv\intl_{\{q\in\pk:\;vqv'<I_m\}} |I_m-vqv
'|^{(\lam-k)/2-d}|q|^{(m-k-1)/2}dq \nonumber \\
&=&2^{-k}\sig_{m,k} \intl_0^{I_k}
|I_k-q|^{(\lam-m-k-1)/2}|q|^{(m-k-1)/2}dq  \nonumber \\
&=&2^{-k}\sigma_{m,k}B_k\left(\frac{m}{2},\frac{\lam-m}{2}\right),
\nonumber \eea and we get the same.

\end{document}